\documentclass{compositio}

\def\b{\mathbb }
\def\phi{\varphi }

\usepackage{amsfonts}
\usepackage{amsmath}
\usepackage{amssymb}
\usepackage{mathrsfs}
\theoremstyle{plain}
\newtheorem{theorem}{Theorem}[section]
\newtheorem{corollary}[theorem]{Corollary}
\newtheorem{lemma}[theorem]{Lemma}
\newtheorem{proposition}[theorem]{Proposition}
\theoremstyle{definition}
\newtheorem{definition}[theorem]{Definition}
\theoremstyle{remark}
\newtheorem*{remark}{Remark}
\newtheorem*{remarks}{Remarks}

\newtheorem{example}[theorem]{Example}
\numberwithin{equation}{section}

\begin{document}

\title{Bessel convolutions on matrix cones}
\author{Margit R\"osler}
\email{mroesler@science.uva.nl}
\address{Ruhr-Universit\"at Bochum\\
Fakult\"at f\"ur Mathematik\\
NA4/30\\
D-44780 Bochum, Germany}

\classification{33C67, 33C80, 43A62, 43A85}
\keywords{Bessel functions, symmetric cones, Dunkl operators, hypergroups, convolution algebras}
\thanks{The author was supported by NWO under project number  B 61-544, and 
by the DFG within SFB/TR 12.}

\begin{abstract}
In this paper we introduce probability-preserving convolution algebras  on cones of positive semidefinite
matrices over one of the division algebras $\b F = \b R, \b C$ or $\b H$ which  interpolate
the convolution algebras of  radial bounded Borel measures on a matrix space $M_{p,q}(\b F)$ 
with $p\geq q$. Radiality in this context means 
invariance under the action of the unitary group $U_p(\b F)$ from the left. 
We obtain a continuous series of commutative hypergroups whose  characters are given by Bessel functions of matrix argument. 
Our  results generalize well-known structures in the rank one case, namely the Bessel-Kingman hypergroups on the positive real line,  to a higher rank setting.
In a second part of the paper, we study structures depending only on the matrix spectra. Under the mapping $r\mapsto \text{spec}(r)$, the 
convolutions on the underlying matrix cone induce a continuous series of 
hypergroup convolutions 
on a Weyl chamber of type $B_q$. The characters  are now Dunkl-type Bessel functions. These convolution algebras on the Weyl chamber naturally extend the harmonic analysis for Cartan motion groups associated with the Grassmann 
manifolds
$U(p,q)/(U_p\times U_q)$ over $\b F$.

\end{abstract}

\maketitle

\section{Introduction}

We start with the 
basic guiding example which corresponds to the rank-one case of our subsequent constructions. For a
natural number $p\geq 2$, consider the set $M_b^{rad}(\b R^p)$  of regular bounded Borel measures 
on $\b R^p$ which are radial, i.e. invariant under orthogonal
transformations. $M_b^{rad}(\b R^p)$ is a commutative Banach algebra with the usual convolution of measures.
Transferring this structure to $\b R_+ = [0,\infty)$ via the mapping
$x \mapsto |x|= (x_1^2 +\ldots + x_p^2)^{1/2}$, 
one obtains a commutative Banach algebra of Borel measures on $\b R_+$. Calculation in polar coordinates shows that
its convolution  $*_p$ is determined on point measures by
\[ \delta_r *_p \delta_s (f) = c_p\int_0^{\pi} f\bigl(\sqrt{r^2 +s^2 -2rs \cos\theta}\bigr) \sin^{p-2}\theta\, d\theta, \quad r,s \in \b R_+, \, f\in C(\b R_+)\]
with a  normalization constant $c_p>0.$  
Now, the above assignment defines a probability measure $ \delta_r *_p \delta_s$ not only for integer $p$ but for all real $p>1$, and it 
extends uniquely to a bilinear and weakly continuous convolution on 
the space $M_b(\b R_+)$ of regular bounded Borel measures on $\b R_+$ 
which is commutative, associative, and probability-preserving (\cite{Ki}).
The interesting point
 about this family of convolutions is that analytic properties
which are valid for integer indices $p$, due to their origin in radial analysis on $\b R^p$, remain true for general indices where no longer any group structure is present. For example, consider the  normalized Bessel functions 
$\, j_\alpha(z) = \,_0F_1(\alpha +1; -z^2/4)\,$
with index $\alpha = p/2-1$. If $p$ is a natural number, then 
\[ j_\alpha (|x|) = \int_{S^{p-1}} e^{-i(x\vert\xi)} d\sigma(\xi), \quad x\in \b R^p\]
($d\sigma$ denotes the normalized Lebesgue surface measure on the unit sphere $S^{p-1}$), and the homomorphism property of the exponential 
function entails the 
 product formula
\[ j_\alpha(r) j_\alpha(s) = \delta_r *_p \delta_s(j_\alpha).\]
This formula, however, extends to arbitrary indices $p>1$, 
see \cite{W}.
The space $\b R_+$ together with the convolution $*_p$ 
is a prominent example of a commutative hypergroup, called the
Bessel-Kingman hypergroup (\cite{BH}). 
A hypergroup is a locally compact Hausdorff space $X$ together with a 
probability preserving convolution of measures on $X$ 
which  generalizes the measure algebra of a locally compact group; in particular, one requires a unit element and the existence of an involution on $X$ generalizing the group inverse in a suitable way; for  
details we refer to \cite{Je} and  Section \ref{orbithypergroups}  below.  There
is a rich harmonic analysis for commutative hypergroups extending the
analysis on locally compact abelian groups. In particular, there is a Haar
measure, a dual space, and a Fourier transform satisfying a Plancherel theorem.
In our example $(\b R_+, *_p)$  the dual space consists of the Bessel functions
\[ \{\phi_s(r) = \,  j_\alpha (rs), \, s\in \b R_+\}, \quad \alpha = \frac{p}{2}-1,\]
and the hypergroup Fourier transform is given by  a Hankel transform,
\[ \widehat f^{\,p}(s) = \, \frac{2^{-p/2}}{\Gamma(p/2)}\int_{\b R_+} f(r) \phi_s(r) r^{p-1}dr.\]
This extends the fact that  the  Fourier transform of a radial function $F(x) = f(|x|) \in L^1(\b R^p)$ is again radial and given by a Hankel transform of $f$
with integral index $p$. 

\medskip

 In the present paper, we generalize the  Bessel convolutions described 
above to a higher rank setting, where the space of ``radii'' is realized as
 a cone of positive semidefinite matrices. More precisely, we construct convolution algebras on such cones which interpolate radial convolution algebras 
on spaces of non-squared matrices. The setting is 
 as follows: For $q\in \b N$ and a natural number $p\geq q$ consider the space
$M_{p,q} = M_{p,q}(\b F)$  of $p\times q$ matrices over one of
the division algebras $\b F = \b R, \b C$ or the quaternions $\b H$. 
It has the structure of a Euclidean vector space with scalar product
$(x|y) = \mathfrak R \text{tr}(x^*y)$ where $x^* = \overline x^t$.
A function (or measure)  on $M_{p,q}$ is called radial if it is invariant under the action of the unitary group $U_p= U_p(\b F)$ on $M_{p,q}$ by left multiplication,
\begin{equation}\label{leftaction}
U_p\times M_{p,q} \to M_{p,q}\,, \quad (u,x) \mapsto ux.
\end{equation}
Observe that this action is via orthogonal transformations, and that  $x$ and $y$ are contained in the same $U_p$-orbit if and only if $x^*x = y^*y$. 
Thus the  space of $U_p$-orbits is naturally parametrized by the cone $\Pi_q = \Pi_{q}(\b F)$ of positive semidefinite $q\times q$-matrices over $\b F$.
If $q=1$ and $\b F= \b R$, then  $\Pi$ coincides with the nonnegative real line 
$\b R_+$ and radiality is equivalent to rotational invariance in $\b R^p$. 
Generalizing the  classical case, there is a radial harmonic analysis 
on $M_{p,q}$ which is based on polar coordinates  with $\Pi_q$ as radial part and 
the Stiefel manifold
\[ \Sigma_{p,q} = \{ x\in M_{p,q}: x^*x = I_q\}\, \cong U_p/U_{p-q}\]
 as a transversal manifold. 
In \cite{FT}, this is developed to some extent
within the general framework of analysis on symmetric cones. 
Indeed, the open cone 
\[\Omega_q  = \{r\in \Pi_{q}: r \text{ strictly positive definite} \}\]
is a symmetric cone within the space
$\,H_q = \{ x\in M_q(\b F): x = x^*\}\,$
 of Hermitian $q\times q$ matrices over $\b F$ which carries 
a natural Euclidean Jordan algebra structure of rank $q$; see Section \ref{anasymm} for details.

For each integer $q\geq p$, we  interpret  radial analysis on $M_{p,q}$ 
in the concise context of a commutative ``orbit hypergroup'' convolution
on the cone $\Pi_q$ which is  derived from the orbit structure  
w.r.t. the action of $U_p$ on $M_{p,q}$.
Similar to the rank-one case, the characters of this hypergroup, i.e. the multiplicative functions which make up the dual,  are obtained
by taking the means of the characters on  $M_{p,q}$ -- usual exponential functions --
 over the Stiefel manifold. This implies
that they are given in terms of Bessel functions $\mathcal J_\mu$ on the 
cone  $\Pi_q$  with (half) integer
 index $\mu = pd/2$,  where  $d=\text{dim}_{\b R}\b F.$ The Bessel function of 
index  $\mu$ is a hypergeometric series of the form
\[\mathcal J_\mu(x) = \sum_{\lambda\geq 0} 
\frac{(-1)^{|\lambda|}}{(\mu)_\lambda\,|\lambda|!} \cdot Z_\lambda(x), \quad x\in H_q;\]
here  the summation goes over all partitions $\{\lambda\in \b N_0^q: \lambda_1\geq \ldots \geq \lambda_q \}$, $(\mu)_\lambda $ is a generalized Pochhammer symbol, and the $Z_\lambda$ are renormalized versions of the  so-called spherical  polynomials associated with the underlying cone. They are homogeneous of degree $|\lambda|$
 and satisfy $\,(\text{tr} x)^k = \sum_{|\lambda|=k} Z_\lambda(x)$ for all $k  \in \b N_0;\,$
c.f.  Section \ref{anasymm}.
The multiplicativity of the Bessel functions ($\mu = pd/2$) with respect to the
orbit convolution on $\Pi_q$ expresses itself in a positive product formula.
Under the technical condition $p\geq 2q$, this  product formula can be written 
in a way which allows analytic continuation with respect to the index $\mu$.
We thus obtain a positive product formula for all
 Bessel functions with index  $\mu  \geq d(q-1/2),$
and an associated continuous series of
hypergroup structures $*_\mu$  on the cone $\Pi_q$ whose dual 
 is given by the functions
\[ \big\{ r\mapsto \phi_s(r) =\, \mathcal J_\mu (\frac{1}{4}rs^2r), \, s\in \Pi_q\big\}.\]
Actually, each hypergroup  $(\Pi_q, *_\mu)$ is  self-dual in a natural way; the neutral element is $0$ and the involution is the identity mapping.
For matrix cones, the Hankel transform of \cite{FT} can now be identified with the
$L^2$-Fourier transform on the underlying hypergroup;
but in addition to the results of \cite{FT}, it is also (and primarily) defined
as a Fourier transform on an $L^1$-convolution algebra. 

\medskip

Before continuing with structural aspects, let us spend some words on 
Bessel functions of matrix argument. They trace back to ideas of Bochner and the fundamental work of Herz \cite{H} and Constantine \cite{Co} in the real case. Much of the interest in these functions is
motivated by questions in  number theory and multivariate statistics. In particular, they occur naturally  in relation with non-central Wishart distributions,
which generalize non-central $\chi^2$-distributions to the higher
rank case, see \cite{Co} and \cite{Mu}.  Nowadays, Bessel functions of matrix argument 
are imbedded into rich theories of multivariable special functions.
First, they can be considered  the $_0F_1$ class among  general $_pF_q$-hypergeometric functions of matrix argument,  where  hypergeometric series are defined in terms of the spherical polynomials; see \cite{GR1} as well as \cite{Ja} for an introduction.
Second, all this can be done in the  general setting of abstract Jordan algebras and symmetric cones, see \cite{FK}. In any case, the  spherical
polynomials depend only on the eigenvalues of their argument. Considered
as functions of the spectra, they can be identified with Jack polynomials of a certain index depending on 
the underlying cone; this was first observed by Macdonald \cite{M2}.
There is a natural theory of hypergeometric expansions 
in terms of Jack polynomials (see \cite{Ka}) which encompasses
the theory on symmetric cones. Finally,
hypergeometric expansions of such kind are intimately related to 
 the modern theory of hypergeometric functions associated with root systems
as developed by Heckman, Opdam, Dunkl and others.
As functions of the spectra, Gaussian hypergeometric functions on a symmetric cone 
 can be identified with
hypergeometric functions associated with a root system of type $BC$ with
a specific choice of parameters
(\cite{BO}). 
Similarly, Bessel functions on a symmetric cone can be considered as a subclass of 
the Bessel functions associated with reduced root systems of type $B$ in the sense of \cite{O},
which play a fundamental role in the theory of rational Dunkl operators \cite{D1}, \cite{D2}.
 This connection 
goes essentially back to \cite{BF} and is made precise in Section \ref{connect} below.

\medskip

In the second part of this paper, we consider 
 structures which depend  only on the spectra
of the matrices from the underlying cone   $\Pi_q$. This amounts to assume invariance under the action of the unitary group $U_q= U_q(\b F)$  by conjugation, 
\[ r\mapsto uru^{-1}, \quad u\in U_q.\]
The orbits under this action are naturally parametrized by the set $\Xi_q$ of possible spectra of matrices from $\Pi_q$,  the eigenvalues being ordered by size:
\[ \Xi_q =\{ \xi= (\xi_1,\ldots, \xi_q) \in \b R^q: \xi_1 \geq \ldots \geq \xi_q\geq 0\}.\]
$\Xi_q$ is a closed Weyl chamber for the reflection group $B_q$
which acts on $\b R^q$ by permutations and sign changes of the coordinates.
Via the canonical mapping from $\,\Pi_q$ onto $\Xi_q$ which assigns to each matrix 
its spectrum, 
the continuous series of hypergroup structures $(\Pi_q, *_\mu)$ with $\mu \geq d(q-1/2)$  induces a series of commutative hypergroup structures $\circ_\mu$ on the chamber $\Xi_q$. The transfer is established by means of so-called orbital mappings. The convolution, Haar measure and  dual space of each  hypergroup on the chamber are made explicit. In particular, the Haar measure of $(\Xi_q, \circ_\mu)$ 
is (up to a constant factor) given by $h_\mu(\xi)d\xi$ where
\[ h_\mu(\xi) = \prod_{i=1}^q \xi_i^{2\gamma +1} \prod_{i <j} (\xi_i^2 - \xi_j^2)^d, \quad \gamma = \mu - \frac{d}{2}(q-1)-1.\]
The hypergroup characters turn out to be certain Dunkl-type Bessel functions 
\[\xi\mapsto J_k^B(\xi, i\eta), \quad \eta\in \Xi_q\]
  associated
with the $B_q$-root system $\{\pm e_i, \,\pm e_i\pm e_j\} \subset \b R^q$.
Here $k$ is a
parameter on the root system which is constant on each subset of
 roots corresponding to a conjugacy class of reflections; in our situation it is
given by $k=(k_1, k_2)$ with 
\[
k_1= \mu - \frac{d}{2}(q-1)-\frac{1}{2} \quad \text {on }\, \pm e_i; \quad
 k_2= \frac{d}{2} \quad \text {on }\, \pm e_i\pm e_j.
\]
The hypergroup  convolution on the Weyl chamber matches the generalized Dunkl translation (\cite{R1}) 
for  Weyl group invariant functions, and we have an interpretation of the  Dunkl transform as 
a hypergroup Fourier transform.

In the geometric cases 
 $\mu = pd/2$,  the support of the probability measure 
$\delta_\xi \circ_\mu\delta_\eta$  on $\Xi_q$ 
describes the set of possible singular spectra of sums $x+y$ with matrices $x,y\in M_{p,q}$ having given singular spectra $\xi$ and $\eta$.
   Further,  the characters are just the bounded spherical functions of the Euclidean type Riemannian symmetric space
$(U_p\times U_q)\ltimes M_{p,q}/(U_p\times U_q)$ associated with the 
Grassmann manifold $U(p,q)/U_p\times U_q$. For general $\mu$
they are characterized, within the theory of rational Dunkl operators,
as the unique analytic solution of a so-called Bessel system; see \cite{O}. It is conjectured
that for arbitrary root systems and non-negative multiplicities, the 
associated Dunkl-type Bessel functions satisfy a positive product formula
and can be characterized as the characters of a commutative hypergroup structure on the 
underlying Weyl chamber. The three continuous series 
($d=1,2,4$) for $B_q$ obtained in this paper are, to our knowledge, 
 the first affirmative examples 
beyond the group cases associated with Cartan motion groups of reductive symmetric spaces.  Some  background, together with
further partial results, is given in  \cite{R2}.

\medskip
\noindent
The organization of the paper is a follows: In Section 2, background on symmetric cones and Bessel functions on cones, as well as some hypergroup analysis are provided. Section 3 
is devoted to the study of Bessel convolutions on matrix cones.
In  Subsection 3.1., orbit convolutions derived from matrix spaces $M_{p,q}$ are  considered. In 3.3., the  corresponding product formula for the involved Bessel functions is analytically extended with respect to the index,  and in  3.4. the associated
series of hypergroup convolutions on the cone are studied. In particular, their  Haar measure and the dual are determined. In 3.5. we analyse 
an interesting critical index. Section 4 is devoted to the induced convolution algebras on a Weyl chamber of type $B$. They are derived 
in  4.1. from the convolutions on the matrix cones and are   then, in the  two final subsections, put into relation
to rational Dunkl theory.

\section{Preliminaries}

In this introductory chapter we provide some relevant background 
on  symmetric cones, in particular matrix cones,  and about 
Bessel functions on such cones. In the main part of the paper, we shall
introduce orbit convolutions and their ``interpolations'' on  matrix cones 
 within the framework of hypergroup theory. As a preparation,
a short account on the relevant notions and facts is included in the present section. For 
a general background on hypergroups, the reader is referred to the fundamental article
\cite{Je} (where the notion "convo" is being used instead of "hypergroup"), 
or to the monograph \cite{BH}.
An excellent  reference for  analysis on symmetric cones is the book \cite{FK},
for special functions on matrix cones see also \cite{GR1} as well as the
classical papers of Herz, James and Constantine, \cite{H}, \cite{Ja}, \cite{Co}.

\subsection{Analysis on symmetric cones}\label{anasymm}

Let $\b F$ be one of the  division algebras $\b F= \b R, \b C$ or the quaternions $\b H$. We denote by $t\mapsto \overline t$  the usual conjugation in $\b F$ and by $\mathfrak R t = \frac{1}{2}(t + \overline t)$ the real part of $t\in \b F$.
 Consider the set of Hermitian $q\times q$-matrices over $\b F$,
\[H_q = H_q(\b F)= \{x\in M_{q}(\b F): x= x^*\}; \quad x^* = \overline x^t.\]
We regard $H_q$ as a Euclidean vector space
with scalar product  $(x\vert y) = \mathfrak R \text{tr}(xy)$, where 
 $\text{tr}$ denotes the trace on $M_q(\b F)$. The dimension of $H_q$ over $\b R$ is 
\[ n= q + \frac{d}{2}q(q-1), \quad d= \text{dim}_{\b R}\b F.\]
With the above scalar product and the Jordan product $x\circ y = \frac{1}{2}(xy + yx),$ the  matrix space $H_q$ becomes a Euclidean (equivalently, a formally real) Jordan algebra with unit $I= I_q$, the unit matrix. The rank of $H_q$, i.e. the number of elements of each Jordan frame in $H_q$, is $q$. 


\medskip
The set $\Omega_q = \Omega_q(\b F)$ of positive definite matrices from $H_q$ 
is a symmetric cone. Recall that a symmetric cone $\Omega$ is
  a proper, non-empty convex cone in a finite-dimensional Euclidean 
vector space 
which is self-dual and homogeneous in the sense that its group of  linear automorphisms acts transitively. Let $G$ denote the connected component of 
this group. Then $K= G\cap O(V)$ is a maximal compact subgroup of $G$ and $\Omega\cong G/K$, 
a Riemannian symmetric space. The matrix cones  are
realized as $\Omega_q(\b F) \cong  GL_q(\b F)/U_q(\b F)$. Hereby $GL_q(\b F)$ acts via $r\mapsto grg^*$, which reduces to conjugation when restricted
to the unitary subgroup $U_q(\b F)$. 
Our main interest will be in the closure of  $\Omega_q$  relative to $H_q$  which coincides with the set of positive semidefinite matrices over $\b F$,
\[ \Pi_q = \Pi_q(\b F)=\,\{x^*x: x\in H_q\} \,=\, \{x^2: x\in H_q\}.\]
 For a general Euclidean Jordan algebra $V,$ the interior $\Omega$ of the set
$\{x^2: x\in V\}$ is a symmetric cone, and  each symmetric cone in a finite-dimensional Euclidean vector space $V$ can be realized in such a way. For details see Section III.3 of \cite{FK}. 
The simple Euclidean Jordan algebras correspond to the irreducible symmetric cones and  are classified. Up to isomorphism,  there are  the above series 
$H_q(\b F)$ with
$\b F= \b R, \b C, \b H$, the exceptional Jordan algebra $H_3(\b O)$, 
as well as  one  infinite series of rank-$2$ algebras corresponding to 
 the Lorentz cones 
\[ \Lambda_n = \{(x^\prime, x_n) \in \b R\times \b R^{n-1}: x_n^2 - 
|x^\prime|^2 >0, \, x_n >0\}.\]

Let $V$ be a simple Euclidean Jordan algebra of rank $q$. 
 We recall that for each $x\in V$ there is a Jordan frame $\{e_1, \ldots, e_q\}$ such that $\, x= \sum_{i=1}^q \xi_ie_i\,$ with real numbers $\xi_i$.
Up to ordering, the $\xi_i$  are uniquely determined  and are called the eigenvalues of $x$. Notice that $x\in \Omega$ iff all its eigenvalues are positive. 
The trace and determinant of $x$ are defined by 
\[\text{tr} \,x = \sum_{i=1}^q \xi_i\,,\quad 
\Delta(x) = \prod_{i=1}^q\xi_i\,.\]
 In the Jordan algebras 
 $H_q(\b F)$, the function $\Delta$  coincides with 
 the usual determinant $\text{det}$
if $\b F= \b R$ or $\b C$ while for $\b F = \b H$ we have $\Delta(x) = (\text{det}\,x) ^{1/2}\, $ when $x$ is being
considered as a complex $2q\times 2q$-matrix in the standard way.

The basic functions for the harmonic analysis on a symmetric cone, and at the same time the basic constituents  for all kinds of 
hypergeometric functions on it, are the so-called spherical polynomials. 
  To write them down, we have to introduce
some more notation:

Consider a fixed Jordan frame  $\{e_1,\ldots, e_q\}$ of $V$. For  $1\leq j\leq q$ and $x\in V$ we denote by $\Delta_j(x)$ the principal minors of $\Delta(x)$ with respect to this frame; in particular $\Delta_q = \Delta$ and each
$\Delta_j$ is a polynomial function on $V$ which is positive on the cone $\Omega$. Recall that a $q$-tuple $\lambda = (\lambda_1, \ldots, \lambda_q)\in \b N_0^q$ 
is called a partition if $\lambda_1\geq \lambda_2\geq\ldots\geq \lambda_q\geq 0$. 
The weight $|\lambda|$ of $\lambda$ is defined as $|\lambda|= \lambda_1 + \ldots + \lambda_q$.
Following a standard convention, we  write
 $\lambda\geq 0$ in order to indicate that $\lambda$ is a partition. For a partition $\lambda$,
 the power function $\Delta_\lambda$ on $V$ is defined by
\[ \Delta_\lambda(x) = \Delta_1(x)^{\lambda_1-\lambda_2} \Delta_2(x)^{\lambda_2-\lambda_3} \cdot\ldots\cdot \Delta_q(x)^{\lambda_q}.\]
It is a homogeneous polynomial of degree $|\lambda|$ and positive on $\Omega$.
In particular, if $x= \sum_{i=1}^q \xi_i e_i$ with $\xi_1, \ldots, \xi_q \in \b R$ then $\Delta_j(x) = \xi_1\cdot\ldots \cdot \xi_j$ and
$\, \Delta_\lambda(x) = \xi_1^{\lambda_1} \cdot\ldots\cdot 
\xi_q^{\lambda_q}.$
 The spherical polynomials of $\Omega$ are indexed by partitions and defined by
\[ \Phi_\lambda(x) = \int_K \Delta_\lambda(kx) dk, \quad \lambda\geq 0,\, x\in V\]
where $dk$ is the normalized Haar measure on $K$. Notice that
$\Phi_\lambda(e) = 1$.
$\Phi_\lambda$ is $K$-invariant and therefore  depends only on the eigenvalues of its argument. The $\Phi_\lambda$  are just the 
polynomial spherical functions of $\Omega\cong G/K$. For the matrix cones over $\b R$  they are known as zonal polynomials, for those over $\b C$ they coincide with the Schur polynomials. Actually, for each 
symmetric cone the 
spherical polynomials are
given in terms of
Jack polynomials (\cite{Sta}). This has been first observed by Macdonald \cite{M2} (see also \cite{F} and the notes in \cite{FK}, Chapt. XI). More precisely, 
let us consider the renormalized polynomials
\[ Z_\lambda(x) := \, d_\lambda \frac{|\lambda|!}{\bigl(\frac{n}{q}\bigr)_\lambda}\cdot \Phi_\lambda(x)\]
where 
$d_\lambda$  denotes the dimension of the vector space of polynomials on $V^{\b C}$ 
which is generated by the elements $z\mapsto \Delta_\lambda(g^{-1}z), \, g\in G.$ (c.f. Sect.XI.5. of \cite{FK}). Let further $C_\lambda^\alpha$ denote the Jack 
polynomials of index $\alpha >0$ (\cite{Sta}),  normalized such that 
\begin{equation}\label{jacksum}
(\xi_1 + \ldots + \xi_q)^k \,=\, \sum_{\lambda\geq 0,\, |\lambda|=k}
 C_\lambda^\alpha(\xi) \qquad \forall \,k\in \b N_0\end{equation}
 (c.f. \cite{Ka}).
Then for $x\in V$ with eigenvalues  $\xi = 
(\xi_1,\ldots ,\xi_q)$ we have 
\[Z_\lambda(x) = C_\lambda^\alpha(\xi)\]
where $\alpha = \frac{2}{d}$ with 
$d= \text{dim}_{\b R}\b F$ for the matrix cones over $\b F$ and $d= n-2$ 
for the Lorentz cone $\Lambda_n$.

In normalization constants, we shall also need the gamma function of a symmetric cone, 
\[ \Gamma_\Omega(z) = \int_{\Omega} e^{-\text{tr} x} \Delta(x)^{z-n/q}dx\]
which converges absolutely for $z\in \b C$ with $\mathfrak R z > \frac{d}{2}(q-1) = \frac{n}{q}-1$ and can be written 
 in terms of the classical gamma function as 
\[ \Gamma_\Omega(z) = (2\pi)^{(n-q)/2}\prod_{j=1}^q\Gamma \Bigl(z - \frac{d}{2}(j-1) \Bigr)\]
see \cite{FK}, Chapt. VII.1.

\subsection{Bessel functions on a symmetric cone}\label{BFM}

Hypergeometric expansions in terms of spherical polynomials have a long history in
multivariate statistics, tracing back to the work of Herz \cite{H}, James \cite{Ja} and Constantine \cite{Co}.
They are important in the study of  Wishart distributions and for questions 
related with total
positivity (\cite{GR2}). In view of the above connection between spherical polynomials and Jack polynomials it is natural to treat 
these classes of functions in the more general framework of
multivariable hypergeometric functions based on Jack polynomial expansions (see \cite{Ka}).
In our context, only hypergeometric functions of type $_0F_1$ (which
are Bessel functions) will be relevant. 

Let $\alpha >0$ be a fixed parameter. For partitions $\lambda=(\lambda_1, \ldots \lambda_q)$ we introduce the generalized Pochhammer symbol 
\[ (c)_\lambda^\alpha =\, \prod_{j=1}^q \bigl(c-\frac{1}{\alpha}(j-1)\bigr)_{\lambda_j}, \quad c\in \b C.\]
For an index $\mu\in \b C$ satisfying $(\mu)_\lambda^\alpha \not= 0$ for all $\lambda\geq 0$ (in particular for $\mu$ with $\mathfrak R \mu >\frac{d}{2}(q-1)$), the generalized hypergeometric function $_0F_1^\alpha(\mu;\,.\,)$ on $\b C^q$ is defined by
\[ _0F_1^\alpha(\mu;\xi) \,=\, \sum_{\lambda\geq 0} 
\frac{1}{(\mu)_\lambda^\alpha |\lambda|!}\cdot C_\lambda^\alpha(\xi).\]
It is known (\cite{Ka}) that this series converges absolutely for all $\xi\in \b C^q$. Similarly, a $_0F_1$-hypergeometric function of two arguments is defined by
\[_0F_1^\alpha(\mu;\xi,\eta) \,=\, \sum_{\lambda\geq 0} 
\frac{1}{(\mu)_\lambda^\alpha |\lambda|!}\cdot \frac{
C_\lambda^\alpha(\xi)C_\lambda^\alpha(\eta)}{C_\lambda^\alpha(\bf 1)}, \quad 
{\bf 1} = (1,\ldots, 1).\] 
Now suppose that $\alpha = 2/d$ where $d$ is the dimension constant of a
simple Euclidean Jordan algebra $V$ of rank $q$ corresponding to the symmetric cone $\Omega$.
Then the Bessel function
$_0F_1^\alpha(\mu;\,.\,)$ essentially coincides with the Bessel function $\mathcal J_\mu$ associated
with  $\Omega$ in the sense of \cite{FK}. Indeed, the latter is defined
by 
\begin{equation}\label{Besselseries}
\mathcal J_\mu(x) =\,\sum_{\lambda\geq 0} 
\frac{(-1)^{|\lambda|}}{(\mu)_\lambda^\alpha |\lambda|!}\,
Z_\lambda(x) \quad \text{for }\, x\in V; \quad \alpha = \frac{2}{d}.
\end{equation}
Thus for $x\in V$ with eigenvalues $\xi = (\xi_1, \ldots, \xi_q)$  we have
\begin{equation}\label{besselid}
\mathcal J_\mu(x) = \,_0F_1^{2/d}(\mu; -\xi).
\end{equation}

\noindent
In Section \ref{spectra} we shall also work with Bessel functions of two arguments $x,y\in V$, 
\[ \mathcal J_\mu(x,y) := \,  \sum_{\lambda\geq 0} \frac{(-1)^{|\lambda|}}{(\mu)_\lambda^\alpha |\lambda|!}
\frac{Z_\lambda(x)Z_\lambda(y)}{Z_\lambda(e)}\]
where $e$ is the unit of $V$. For  $x,y\in V$ with eigenvalues $\xi = (\xi_1, \ldots ,\xi_q)$ and $\eta = (\eta_1, \ldots, \eta_q)$ respectively, we thus have
\begin{equation}\label{dunklbessel2} \mathcal J_\mu(x,y) = \,_0F_1^{2/d}(\mu;i\xi,i\eta).
\end{equation}

\noindent
The following estimate is especially useful for large indices $\mu$.

\begin{lemma}\label{F_estimate} For $\alpha >0$ let  $_0F_1^\alpha(\mu;\,.\,)$ denote 
the associated generalized Bessel function of index $\mu$ on $\b C^q$.
Suppose that $\mathfrak R \mu > \frac{1}{\alpha}(q-1)$. Then for
$\xi=(\xi_1,\ldots,\xi_q)\in \b C^q$, 
\[ \big|\phantom{}_0F_1^\alpha(\mu;\xi)\big| \,\leq \, e^{|\xi_1|+\ldots +|\xi_q|}.\]
\end{lemma}

\begin{proof}
For $n\in \b N_0$ and $z\in \b C$ with $\mathfrak R z >0$ we have
\[ \Big| \frac{\Gamma(z)}{\Gamma(z+n)}\Big| \,\leq\, \frac{\Gamma(\mathfrak R z)}{\Gamma(\mathfrak R z + n)}\,\leq \, 1.\]
Therefore
$\, \big| (\mu)_\lambda^\alpha\big| \geq 
\big| (\mathfrak R\mu)_\lambda^\alpha\big|\geq 1.\,$
Moreover, it is known that the coefficients of the Jack polynomial $C_\lambda^\alpha$ 
in its monomial expansion are all non-negative (see \cite{KnS}). This implies that
\[ \big| C_\lambda^\alpha(\xi)\big| \,\leq C_\lambda^\alpha(|\xi_1|,\ldots, |\xi_q|).\]
Using relation \eqref{jacksum}, we therefore obtain
\[ \big|\phantom{}_0F_1^\alpha(\mu;\xi)\big| \,\leq \sum_{\lambda\geq 0} \frac{1}{|\lambda|!} C_\lambda^\alpha(|\xi_1|,\ldots, |\xi_q|)\,=\, e^{|\xi_1|+\ldots +|\xi_q|}.\]
\end{proof}

\begin{corollary} \label{B_estimate}
Let $\mathcal J_\mu$ denote the Bessel function of index $\mu$ 
associated with the symmetric cone $\Omega$ inside the Jordan algebra $V$ of rank $q$.
Suppose that $\mathfrak R \mu > \frac{d}{2}(q-1)$. Then for $x\in V$ with eigenvalues $\xi_1, \ldots \xi_q$, 
\[ \big|\mathcal J_\mu(x) \big|\,\leq \, e^{|\xi_1| +\ldots + |\xi_q|}.\]
\end{corollary}

We mention at this point that for $\mathfrak R\mu > d(q-1) +1$ and $a\in \Omega$,
the Bessel function on $\Omega$ has the absolutely convergent integral representation (an inverse Laplace integral)
\[ \mathcal J_\mu(r) = \frac{\Gamma_\Omega(\mu)}{(2\pi i)^n} \int_{a + iV} 
e^{\text{tr} w} e^{-(w^{-1}\vert r)} \Delta(w)^{-\mu} dw \quad (r\in \Omega)\]
see Prop. XV.2.2. of \cite{FK}. This integral was originally used by Herz
\cite{H} to define Bessel functions on the cones of positive definite matrices over $\b R$; see also \cite{FT} and \cite{Di} for arbitrary symmetric cones. 
Corollary \ref{B_estimate} slightly improves the order estimate for $\mathcal J_\mu$ 
given in \cite{H}, p.486 for $\b F = \b R$. 


\subsection{Hypergroups}\label{orbithypergroups}

We start with some notation:
For a locally compact Hausdorff space  $X$, let  $M_b(X)$ denote the Banach space of all bounded
regular (complex) Borel measures on $X$ with total variation norm, and $M^1(X)\subset M_b(X) $ the set of all probability measures.
With $\delta_x$  we denote the point
measure in $x\in X$.
We use the notions $C(X), C_b(X)$ and $C_c(X)$ for the spaces of continuous complex-valued functions on $X$, those which are bounded,  and those having compact support respectively. Further, $C_0(X)$ is the set of functions from $C(X)$ which vanish at infinity.

\begin{definition} A hypergroup 
$(X,*)$ 
is a locally compact Hausdorff space $X$ with a bilinear and associative  convolution $*$  on $M_b(X)$ with the following properties:
\begin{enumerate}
\item[\rm{(1)}] The map $(\mu,\nu) \mapsto \mu*\nu$ is weakly continuous, i.e. w.r.t. the topology induced by $C_b(X)$.
\item[\rm{(2)}] For all $x,y\in X$, the product $\delta_x*\delta_y$ of point measures
is a compactly supported probability measure on $X$. 
\item[\rm{(3)}] The mapping $(x,y)\mapsto \text{supp}(\delta_x *\delta_y) $ from $X\times X$ into the space of nonempty compact subsets of $X$ is continuous with respect to the Michael topology (see \cite{Je}). 
\item[\rm{(4)}] There is a neutral element $e\in X$, satisfying $\delta_e*\delta_x = \delta_x * \delta_e
= \delta_x$ for all $x\in X$.
\item[\rm{(5)}] There is a continuous involutive automorphism $x\mapsto\overline x$ on $X$ such that $\delta_{\overline x} * \delta_{\overline y} =(\delta_y*\delta_x)^-$ and $\,x=\overline y \Longleftrightarrow e\in \text{supp}(\delta_x*\delta_y).$
(For $\mu\in M_b(X)$, the measure  $\mu^-$ is given by
$\mu^-(A)=\mu(\overline A)$ for Borel sets $A\subseteq X$).
\end{enumerate}
A hypergroup $(X,*)$ is called commutative if its convolution is commutative.
\end{definition}


For a hypergroup $(X,*)$, the space $(M_b(X),*)$ is a 
Banach algebra with unit $\delta_e$. Notice that
by the density of the finitely supported measures in $M_b(X)$, the convolution is uniquely determined as soon as it is given for point measures.

Of course every locally compact group is a hypergroup with the usual group convolution. In our context, only commutative hypergroups will be relevant. 
Of particular interest will be the following

\begin{example}(\cite{Je}, Chapt. 8) 
\label{ex_1}Let $(G,+)$ be a locally compact abelian group
and $K$ a compact subgroup of $\text{Aut}\,G.$ 
Then the space $\,G^K = \{K.x: x\in G\}$ 
 of $K$-orbits in $G$ is a locally compact Hausdorff space with the quotient topology and becomes a commutative hypergroup with the (natural) definition
\[ (\delta_{K.x} * \delta_{K.y})(f) =\,\int_K f(K.(x+ky))dk, \quad f\in C_c(G^K).\]
$(G^K, *)$ is called a orbit hypergroup; its neutral element is $K.0=0$ and the involution is $(K.x)^{-} = K.(-x).$
\end{example}

In the following, we collect some further ingredients 
underlying the harmonic analysis on a commutative hypergroup $X$.
There exists (up to normalization) a unique Haar measure $\omega$ on $X$, i.e. a 
positive Radon measure satisfying
\[ \int_X f(x*y)d\omega(y) = \int_X f(y)d\omega(y) \quad\forall\, x\in X, f\in C_c(X)\]
where we use the notation
\[f(x*y):= \int_X f\> d(\delta_x*\delta_y).\]
Similar to the dual  of a locally compact abelian group, 
one defines the dual space of $X$ as 
\[\widehat X:=\{\phi\in C_b(X):\> \phi\ne 0,\>\> \phi(\bar x) = \overline{\phi(x)} 
\>\text{ and }\,  \phi(x*y)=
 \phi(x)\phi(y) \, \forall \, x,y\in X\}.\]
 This is a locally compact Hausdorff space with the topology of 
  compact-uniform convergence.  
The elements of  $\widehat X$ are also called the
characters of $X$. They are the constituting functions for harmonic analysis
on $(X,*)$. The  Fourier transform on $M_b(X)$ is defined by
\[\widehat \mu(\phi):= \int_X \overline{\phi(x)}\>
d\mu(x), \quad\phi\in \widehat X,\]
and on $ L^1(X, \omega)$ by $\widehat f:= \widehat{f\omega}$.
The Fourier transform is injective, and there exists a 
unique positive Radon measure $\pi$ on $\widehat X$, called the Plancherel measure of $(X,*)$, such that $f\mapsto \widehat f$ extends to an isometric 
isomorphism from $L^2(X,\omega)$ onto $L^2(\widehat X,\pi)$. 
While the Haar measure $\omega$ has full support, the support of $\pi$ 
may be a proper subset of $\widehat X$ only. 

\begin{example}\label{ex_2} [Continuation of example \ref{ex_1}.]
For a Haar measure $m$ on the group $G$, the image measure of $m$ under the the canonical map
$\pi: G\to G^K$ provides a Haar measure on the orbit hypergroup $G^K$.
Further, it is easily seen that 
the functions
 \[ \phi_\alpha(K.x) :=
 \int_K \alpha(k.x) dk, \quad \alpha\in \widehat G,\]
 belong to the dual of $G^K$. 
 Actually, we have
 \end{example}
 
 \begin{lemma}\label{dual}
\begin{enumerate} 
\item[\rm{1.}] $\displaystyle\widehat{G^K} = \{ \phi_\alpha : \alpha\in \widehat G\}.$
\item[\rm{2.}] $\phi_\alpha = \phi_{\alpha^\prime}$ iff $\alpha$ and $\alpha^\prime$ are contained in the same 
orbit under the dual action of $K$ on $\widehat G$, given by $(k.\alpha)(x) = \alpha(k^{-1}.x)$
\end{enumerate}
\end{lemma}

\begin{proof}
It is known that $\widehat{G^K}$  can be identified with
the extremal points of the set
\[ \Xi = \{ \beta\in C(G):\, \beta \text{ positive definite and $K$-invariant with }\, \beta(0) = 1\};\]
the identification being given by $\,\beta \mapsto \phi_\beta, \, \phi_\beta(K.x) = \beta(x)$ for $x\in G$,
see \cite{Ro}.
By Bochner's theorem, each $\beta\in \Xi$ is of the form
$\beta(x) = \int_{\widehat G} \alpha(x)d\mu(\alpha)\,$ with some $K$-invariant 
probability measure $\mu$ on $\widehat G$. In this way, the extremal points of $\Xi$ 
correspond to the measures of the form $\mu_\alpha =\int_K \delta_{k.\alpha} dk\,$ with $\alpha\in \widehat G$. Further, $\mu_\alpha = \mu_{\alpha^\prime}$ iff $\alpha$ and $\alpha^\prime$ belong to the same $K$-orbit in $\widehat G$. 
\end{proof}

A commutative hypergroup $(X,*)$ is called self-dual
if there exists a homeomorphism $\Psi: X\to \widehat X$ such that
\[ \Psi(x)(z)\Psi(y)(z)\,=\,\int_X\Psi(w)(z)d(\delta_x*\delta_y)(w) \quad\forall 
x,y,z\in X.\]
In that case, $\widehat X$ is a dual hypergroup in the sense
of \cite{Je}, 12.4. where the convolution product
$\delta_{\Psi(x)}* \delta_{\Psi(y)}$ is just the image measure of $\delta_x*\delta_y$ under $\Psi$. It is easily checked that the image measure $\Psi(\omega)$ of the Haar measure $\omega$ on $X$  is a Haar measure on $\widehat X$. Thus by 
Theorem 12.4A of \cite{Je}, $\Psi(\omega)$ coincides up to a multiplicative constant with the Plancherel measure $\pi$ of $X$.

\section{Convolution structures associated with Bessel functions on a matrix cone}
\label{convocones}

\subsection{Orbit hypergroups on a matrix cone}\label{orbitcone}

For natural numbers $p\geq q$, consider  the matrix space $M_{p,q}=M_{p,q}(\b F)$ of 
$p\times q$-matrices over $\b F$. We regard $M_{p,q}$ as a real vector space,  equipped with the Euclidean scalar product
 $(x\vert y) := \mathfrak R\text{tr}(x^*y)$ and norm $\,\|x\| = \sqrt{\text{tr}(x^*x)}$.  Here $\mathfrak R t = \frac{1}{2}(t+ \overline t)$ is the real part of $t\in \b F$  and 
$\text{tr}$ denotes the trace in $M_q(\b F)$. In the square case $p=q$, $\|.\|$ is just the Frobenius norm.  Let us consider the action of the
 unitary group $U_p$  on $M_{p,q}$ by left multiplication,  
\[ U_p\times M_{p,q} \to M_{p,q}\,, \quad (u,x) \mapsto ux.\] 
The
 orbit space 
$M_{p,q}^{U_p}$ for this action
can be identified with the space 
$\Pi_q= \Pi_q(\b F)$ of positive semidefinite $q\times q$ matrices over $\b F$ via 
\[ U_p.x \mapsto \,\sqrt{x^*x}\, =:\, |x|.\]
 Here for $r\in \Pi_q\,,\, \sqrt{r}$ denotes the unique positive semidefinite square root of $r$.  It is easy to see
that the above bijection becomes a homeomorphism 
 when $M_{p,q}^{U_p} $
 is equipped with the quotient topology and $\Pi_q$
with the subspace topology induced from $M_q$. 
(Indeed, both the canonical map $x\mapsto U_p.x$ from
$M_{p,q}$ onto $M_{p,q}^{U_p}$ and the mapping $x\mapsto \sqrt{x^*x}, \, M_{p,q}\to \Pi_q\,$ are open and continuous.) Notice that  the Stiefel manifold
\[ \Sigma_{p,q} = \{x\in M_{p,q}:\, x^*x = I_q\}\]
is the orbit of the block matrix
\[ \sigma_0 :=  \begin{pmatrix} I_q\\ 0
\end{pmatrix} \in M_{p,q}.\]

Before calculating the orbit hypergroup convolution for the above action,
 let us recall
the basic aspects of radial analysis in $M_{p,q}$ as developed in \cite{FT}.
As in the introduction, radiality here means invariance under the above action
of $U_p$; a function $F$ on $M_{p,q}$ is radial iff it is of the form $F(x) = f(|x|)$ for some $f:\Pi_q\to \b C$.
Suitable  polar coordinates in  $M_{p,q}$ are defined as follows:
Let  $d\sigma$ denote the $U_p$-invariant measure on $\Sigma_{p,q}$, normalized according to $\int_{\Sigma_{p,q}}d\sigma = 1,$ and let
\[ M_{p,q}^\prime = \{x\in M_{p,q}:\, \Delta(x^*x)\not= 0\}\]
which is open and dense in $M_{p,q}.$ Then  the mapping
\[ \Omega_q \times \Sigma_{p,q} \to M_{p,q}^\prime\,, \quad (r,\sigma) \mapsto \sigma\sqrt{r}\]
is a diffeomorphism, and for integrable functions $f:M_{p,q} \to \b C$ one has
\[ \int_{M_{p,q}} f(x)dx \,=\, C_{p,q} \int_{\Omega_q}\int_{\Sigma_{p,q}} f(\sigma\sqrt{r}) \Delta(r)^\gamma dr\, d\sigma\]
with 
\begin{equation}\label{const1}C_{p,q} = \frac{\pi^{dpq/2}}{\Gamma_{\Omega_q}\bigl(\frac{dp}{2}\bigr)}\, \quad\text{and} \quad
\gamma = \frac{dp}{2} - \frac{n}{q}.\end{equation} As before, $n$ is the real dimension of $H_q$.
Our notion slightly differs from that of \cite{FT} 
(and the monograph \cite{FK}); it is adapted to the left action \eqref{leftaction} of  $U_p$ on $M_{p,q}$ while the notion of \cite{FT} would require to have $U_p$ acting from the right on $M_{q,p}$. 
 As Lebesgue measure on $M_{p,q}$ is $U_p$-invariant, the Fourier transform of a radial function in
$ L^1(M_{p,q})$ is again radial, and also the convolution of two radial
functions is radial; both can be calculated by use of polar coordinates, c.f. \cite{FT}. We shall come back to this shortly.  Following Example 
\ref{ex_1}, we obtain the orbit hypergroup convolution 
\[ (\delta_r * \delta_s)(f) = \int_{U_p} 
f\bigl(\vert \sigma_0 r + u \sigma_0 s\vert\bigr) du, \quad r,s\in \Pi_q.\]
The image measure of $du$ under the mapping $\, U_p\to \Sigma_{p,q}\,,\,  u\mapsto u\sigma_0$ is  $U_p$-invariant and therefore coincides with the $U_p$-invariant measure $d\sigma$. Hence,
\begin{equation}\label{convo}
(\delta_r * \delta_s)(f) = \int_{\Sigma_{p,q}} f(\vert\sigma_0 r+ \sigma s\vert ) d\sigma\,
= \,\int_{\Sigma_{p,q}} f\bigl(\sqrt{r^2 + s^2 + r\widetilde\sigma s +(s\widetilde\sigma r)^*}\bigr) d\sigma
\end{equation}
where $\,\widetilde\sigma = \sigma_0^*\sigma$ is the $q\times q$-matrix whose rows are given by the first $q$ rows of $\sigma$. 
Actually, this convolution depends on $p$ (and $q$), which we 
surpress for the moment.
 The neutral element of the orbit hypergroup $(\Pi_q, *)$ is $0$, and the involution is the identity mapping (because $x\in M_{p,q}$ and $-x$ are in the same $U_p$-orbit).
Further, according to  Lemma \ref{dual},
the dual space of $(\Pi_q, *)$ consists of the functions $\phi_s, \,s\in \Pi_q$ with  
\begin{equation}\label{Bochnerrep}
 \phi_s(r) = \,\int_{U_p} e^{-i(u\sigma_0 r\vert \sigma_0 s)} du\, =\,
 \int_{\Sigma_{p,q}} e^{-i(\sigma\vert \sigma_0 sr)} d\sigma .
\end{equation}
These are Bessel functions. Indeed,  
according to Propos. XVI.2.3. of \cite{FK} we have for $x\in M_{p,q}$ the 
identity
\begin{equation}\label{intrepfaraut}\int_{\Sigma_{p,q}} e^{-i(\sigma\vert x)}d\sigma = \mathcal J_\mu 
\bigl(\frac{1}{4}x^*x\bigr), \quad  \mu = \frac{pd}{2}
\end{equation}
where $\mathcal J_\mu$ is the Bessel function of index $\mu$ associated with the symmetric cone 
$\Omega_q$ as in Section \ref{BFM}. 
Thus
\[ \phi_s(r) = \mathcal J_\mu\bigl(\frac{1}{4}rs^2r).\]
We mention at this point that up to a constant factor, 
$\phi_s(r)$ coincides with the Bessel function
$J(r^2, s^2)$ in \cite{FT}.
As $\mathcal J_\mu(x)$ depends only on the eigenvalues of $x$ and 
as the matrices $rs^2r$ and $sr^2s$ have the same eigenvalues,
we see that
\[ \phi_s(r) = \phi_r(s) \quad \forall r, s\in \Pi_q.\]
This implies that the hypergroup $(\Pi_q, *)$ is self-dual via $\,r\mapsto
\phi_r$. In order to stress the dependence of its convolution on $p$ (or equivalently, on the index 
$\mu= pd/2$), we denote it by $*_\mu$ from now on and write $\Pi_{q,\mu}$ for the  orbit hypergroup $(\Pi_q, *_\mu)$.

A Haar measure $\omega_\mu$ on the hypergroup $\Pi_{q,\mu}$ is obtained by taking the image measure of the (normalized)
Lebesgue measure $(2\pi)^{-pqd/2}dx$ on $M_ {p,q}$ under the mapping $x\mapsto |x|.$ Using again polar coordinates, we obtain
\[ \omega_\mu(f) = \frac{2^{-\mu q}}{\Gamma_{\Omega_q}(\mu)}\int_{\Omega_q} f(\sqrt r)\Delta(r)^\gamma dr, \quad f\in C_c(\Pi_q).\]

Fourier transform and convolution of radial functions on $M_{p,q}$ are calculated in our notion a follows:
Suppose $F,G \in L^1(M_{p,q})$ are radial with $F(x) = f(|x|)$ and $G(x) = g(|x|)$. Then the Fourier transform of $F$ is
\[ \widehat F (\lambda) \,=\, \frac{1}{(2\pi)^{dpq/2}}\int_{M_{p,q}} F(x) e^{-i(\lambda\vert x)} dx \,=\, \int_{\Pi_q} f(r) \phi_{|\lambda|}(r) d\omega_\mu(r)\,=\, \widehat f(|\lambda|).
 \]
The convolution of $F$ and $G$ is given by  $F*G(x) = H(|x|)$ with
 \begin{align} H(r) = &\,(2\pi)^{-dpq/2} \int_{M_{p,q}}F(\sigma_0 r -y)G(y) dy 
\notag\\ =\,&\frac{2^{-\mu q}}{\Gamma_{\Omega_q}(\mu)} \int_{\Omega_q}
\Bigl(\int_{\Sigma_{p,q}} f(|\sigma_0 r -\sigma\sqrt{s}\,|)\Bigr) G(\sqrt{s})
\Delta(s)^\gamma ds \notag\\
=\,&\int_{\Pi_q} \delta_r *_\mu \delta_s(f)\, g(s) d\omega_\mu(s) \,=\, (f*_\mu g)(r).
\notag
\end{align}

The multiplicativity of the characters $\phi_s$ implies a positive product formula for the Bessel functions $\mathcal J_\mu$ with index $\mu = pd/2, \,p\geq q$ an integer:
\begin{equation}\label{prod2} \mathcal J_\mu(r^2)\mathcal J_\mu(s^2) = \int_{\Sigma_{p,q}} 
    \mathcal J_\mu\bigl(r^2 + s^2 +r\widetilde\sigma s + s\widetilde\sigma^*r\bigr)
d\sigma \quad \forall r,s\in \Pi_q.
\end{equation}
We shall generalize  this  formula to Bessel functions of arbitrary index 
$\mu\in \b C$ with $\mathfrak R\mu > d(q-1/2)$. For real indices we shall obtain
a positive
product formula, which leads to  a continuous family 
of hypergroup structures on $\Pi_q$
 beyond those which have 
a realization as orbit hypergroups as above. The decisive observation towards this aim is  that only the reduced matrix $\widetilde\sigma$ 
occurs in the integrands of \eqref{prod2} and \eqref{convo}.
 In the following section, we shall introduce coordinates on the 
Stiefel manifold which are adapted to this situation.

\subsection{Split coordinates on the Stiefel manifold}

Let $k\in \b N$ with $p-k\geq q$. We decompose $\sigma\in \Sigma_{p,q}$ 
as $\sigma= \begin{pmatrix} v\\ w\end{pmatrix}$ with $v\in M_{k,q}$ and $w\in M_{p-k,q}$. For fixed $q$, put
\[ D_{k} = \{ v\in M_{k,q}: v^*v < I\}\]
where the notion $x<y$ for $x,y\in M_q(\b F)$ means that $y-x$ is (strictly) positive-definite. 

\begin{proposition}\label{polar1}
The mapping
\[ \Phi: D_{k}\,\times \,\Sigma_{p-k,q} \,\longrightarrow \, \Sigma_{p,q}, \quad (v, \sigma^\prime)  \, \longmapsto \, \begin{pmatrix} v\\
 \sigma^\prime\sqrt{I-v^*v}\end{pmatrix} \]
is a diffeomorphism onto a dense and open subset $\widetilde\Sigma_{p,q}$ 
of $\Sigma_{p,q}$.
              Let $d\sigma$ and $d\sigma^\prime$ denote the normalized 
Riemannian volume elements on
$\Sigma_{p,q}$ and $\Sigma_{p-k,q}$ respectively, and let
$\eta := \frac{d}{2}(p-k) -\frac{n}{q}$. 
Then on $\widetilde\Sigma_{p,q}$, 
\[ d\sigma = c_{k,p}\cdot \Delta(I-v^*v)^\eta \,d\sigma^\prime dv \]
where
\[ c_{k,p} = \bigl(\int_{D_{k}}\Delta(I-v^*v)^\eta dv\bigr)^{-1} \,=\,\pi^{-dkq/2} \frac{\Gamma_{\Omega_q}\bigl(\frac{dp}{2}\bigr)}{\Gamma_{\Omega_q}\bigl(\frac{d(p-k)}{2}\bigr)}.\]
\end{proposition}

In the important special case $p\geq 2q, k=q$ with $\b F = \b R$, this result
goes back to \cite{H}. For the reader's 
convenience we nevertheless supply a concise proof along a different approach.

\begin{proof} It is easily checked that $\Phi:  D_{k}\,\times\,\Sigma_{p-k,q} \,\rightarrow \, \Sigma_{p,q}$ is $C^1$ and injective, and that its image
$\widetilde\Sigma_{p,q}$ is dense in $\Sigma_{p,q}$.
Let $f\in C(\Sigma_{p,q})$ and extend it to $F\in L^1(M_{p,q})$ by
$\,F(\sigma\sqrt r):= f(\sigma)$ if $\frac{1}{2}I <r< \frac{3}{2}I$ and $F=0$ else.
Then 
\begin{equation}\label{part1} \int_{M_{p,q}} F(x)dx\,  =\, C\int_{\Sigma_{p,q}}f d\sigma \quad\text{ with }\, C = C_{p,q}\cdot\int_{\frac{1}{2}I<r<\frac{3}{2}I} \Delta(r)^\gamma dr. 
\end{equation}
On the other hand, we
write $x\in M_{p,q}$ in block form as 
$x = \begin{pmatrix} x_1\\ x_2\end{pmatrix}$ 
with $x_1\in M_{k,q}$ and $x_2\in M_{p-k,q}$.
Then
\[
 \int_{M_{p,q}}\! F(x)dx \,=\,
C_{p-k,q}\int_{M_{k,q}}\int_{\Omega_q}\int_{\Sigma_{p-k,q}}F\begin{pmatrix}x_1\\ \sigma_2\sqrt{r_2}\end{pmatrix}\,\Delta(r_2)^\eta \, dr_2 d\sigma_2 dx_1
\]
where we used polar coordinates $x_2 = \sigma_2\sqrt{r_2}$ in $M_{p-k,q}$. 
Consider now the successive transforms 
$\, r_2 \mapsto x_1^*x_1 + r_2 =: r\,$ with $x_1$ fixed
and $\, x_1\mapsto v:= x_1 r^{-1/2}\,$ with $r$ fixed.
We have $dr_2 = dr$ and, by Lemma 2 of \cite{FT}, $dx_1 = \Delta(r)^{dk/2}dv.$
Thus in a first step, 
\[
 \int_{M_{p,q}}\! F(x)dx \,=\,
  C_{p-k,q}\int_{M_{k,q}}\int_{x_1^*x_1+ \Omega_q}\int_{\Sigma_{p-k,q}} F\begin{pmatrix}x_1\\ \sigma_2\sqrt{r-x_1^*x_1}\end{pmatrix}\Delta(r-x_1^*x_1)^\eta \, dr d\sigma_2 dx_1\]
where $\, F\begin{pmatrix}x_1\\ \sigma_2\sqrt{r-x_1^*x_1}\end{pmatrix} \,=\, f\begin{pmatrix}v\\ \sigma_2 w\end{pmatrix}$ with $\,w= 
\sqrt{r-\sqrt{r} v^*v \sqrt{r}}\cdot r^{-1/2}$. As $\, w^*w = I -v^*v$, it follows that $w= u\sqrt{I -v^*v}$ 
with some $u\in U_q$. The invariance of the Lebesgue measure on $M_{p-k,q}$ under the action of $U_q$  by right multiplication easily implies that
$d\sigma_2$ is  invariant under this action of $U_q$ on $\Sigma_{p-k,q}.$
In view of the identity  $\eta + kd/2 = \gamma$, the above integral therefore becomes
\begin{align}
  C_{p-k,q}&\int_{\frac{1}{2}I<r<\frac{3}{2}I}\int_{D_{k}}\int_{\Sigma_{p-k,q}} 
    f\begin{pmatrix}v\\ \sigma_2\sqrt{I-v^*v} \end{pmatrix} 
\Delta(r)^{\eta + kd/2} \Delta(I-v^*v)^\eta \, d\sigma_2 dv\, dr\,\notag\\
=\,& C\cdot\frac{C_{p-k,q}}{C_{p,q}}\int_{D_{k}}\int_{\Sigma_{p-k,q}}
  f\begin{pmatrix}v\\ \sigma_2\sqrt{I-v^*v}\end{pmatrix} \Delta(I-v^*v)^\eta \, d\sigma_2 dv.\notag \end{align}
Together with \eqref{part1}, this  gives the stated Jacobian (which does not vanish on $D_{k}\,\times\,\Sigma_{p-k,q}$) and
proves the claimed diffeomorphism property of $\Phi$.
\end{proof}

If $p\geq 2q$, then the above proposition with $k=q$ leads to a nice integration formula for
functions $f(\sigma)$ on $\Sigma_{p,q}$ which depend only on the first $q$ rows
of $\sigma$.
For abbreviation, we put 
\[ \varrho:= d\bigl(q-\frac{1}{2}\bigr) +1\]
and for $\mu\in \b C$ with $\mathfrak R\mu >\rho -1$, 
\[ \kappa_\mu := \int_{D_q} \Delta(I-v^*v)^{\mu-\rho} dv.\]
The explicit value of $\kappa_\mu$ is obtained by using polar coordinates and Thm. VII.1.7. of \cite{FK} about beta integrals on symmetric cones,
\[ \kappa_\mu\, =\,  C_{q,q} \int_{r<I} \Delta(I-r)^{\mu-\varrho}\Delta(r)^{d/2-1}dr\,=\,\pi^{dq^2/2}\cdot\frac{\Gamma_{\Omega_q}\bigl(\mu - \frac{dq}{2}\bigr)}{\Gamma_{\Omega_q}(\mu)}.   \]

\begin{corollary}\label{red} Let $p\geq 2q$. Then for $f\in C(\Sigma_{p,q})$ of the form
\[ f(\sigma) = \widetilde f(\widetilde\sigma), \quad \widetilde\sigma = \sigma_0^*\sigma\]
one has
\[ \int_{\Sigma_{p,q}} fd\sigma \,=\,\frac{1}{\kappa_{pd/2}} \int_{D_q} \widetilde f(v)\,\Delta(I-v^*v)^{pd/2-\varrho} dv.\]
\end{corollary}

\noindent
The dependence on $p$ now occurs only in  the density, not in the domain of integration.

\subsection{A product formula for Bessel functions
with continuous index}

Remember that the integrand in the convolution formula \eqref{convo} and the product formula \eqref{prod2}
depend only on the reduced matrix $\widetilde\sigma = \sigma_0^*\sigma$.
By Corollary \ref{red} we obtain

\begin{proposition}\label{convodisc}
Suppose that $p\geq 2q$ and let $\mu = pd/2$.
\begin{enumerate}
\item[\rm{(1)}] 
  The convolution \eqref{convo} (i.e. $*_\mu$) on $\Pi_q$ can be written as
\[(\delta_r *_\mu \delta_s)(f) = \frac{1}{\kappa_\mu}\int_{D_q} f\bigl(\sqrt{r^2 + s^2 + rvs + sv^*r}\,\bigr)\,
\Delta(I-v^*v)^{\mu-\varrho}\, dv. \]
\item[\rm{(2)}] The Bessel function $\mathcal J_\mu$ satisfies the product formula
\[ \mathcal J_\mu\bigl(r^2\bigr) \mathcal J_\mu\bigl(s^2\bigr) \,=\, 
\frac{1}{\kappa_\mu} \int_{D_q} \mathcal J_\mu\bigl(r^2+s^2 + rvs + sv^*r\bigr)\,\Delta(I-v^*v)^{\mu- \varrho}dv \quad \forall\, r,s\in \Pi_q.\]
\end{enumerate}
\end{proposition}

 We are now going to extend the integral formulas of Proposition \ref{convodisc}
 to arbitrary indices $\mu$ within the half plane $\{\mu\in \b C: \mathfrak R \mu >
  \varrho-1 = d\bigl(q-\frac{1}{2}\bigr)\}.$ We will use the standard technique
  (c.f. \cite{Ste} for the rank one case), namely analytic continuation with respect to $\mu$.  The argumentation will be based on a classical theorem of Carlson:

\begin{theorem}[\cite{T}, p.186] Let $f(z)$ be holomorphic in a neighbourhood of
$\{z\in \b C:\mathfrak R z \geq 0\}$ satisfying $f(z) = O\bigl(e^{c|z|}\bigr)$
on $\,\mathfrak R z \geq 0$ for some $c<\pi$. If $f(z)=0$ for all $z\in \b N_0$, then $f$ is identically zero.
\end{theorem}

\noindent
The following theorem is the main result of this section.

\begin{theorem}\label{main1}
Let $\mu\in\b C$ with $\mathfrak R\mu > \varrho-1 = d\bigl(q-\frac{1}{2}\bigr).$ Then the Bessel function $\mathcal J_\mu$ satisfies the product formula
\begin{equation}\label{prod3}
\mathcal J_\mu\bigl(r^2\bigr) \mathcal J_\mu\bigl(s^2\bigr) \,=\, 
\frac{1}{\kappa_\mu} \int_{D_q} \mathcal J_\mu\bigl(r^2+s^2 + rvs + sv^*r\bigr)\,\Delta(I-v^*v)^{\mu- \varrho}dv \quad \forall r,s\in \Pi_q.
\end{equation}   
\end{theorem}

\begin{remark}
It is easy to obtain from \eqref{prod3}  the more general identity
\begin{equation}\label{extended} \mathcal J_\mu\bigl(x^*x\bigr) \mathcal J_\mu\bigl(y^*y\bigr) \,=\, 
\frac{1}{\kappa_\mu}\int_{D_q} \mathcal J_\mu\bigl(x^*x+ y^*y + x^*vy + y^*v^*x\bigr) \Delta(I-v^*v)^{\mu- \varrho}dv 
\end{equation}
for arbitrary $x,y\in M_q$. For this, recall that every $x\in M_q$ 
has a polar decomposition $x = u|x| $ with
$|x| = \sqrt{x^*x}\in \Pi_q$ and a unitary matrix $u\in U_q$.
Choose $r=|x|=u^*x$ and $s=|y|=w^*y$ with $u,w \in U_q$ in \eqref{prod3}.
Then 
\[ \mathcal J_\mu\bigl(x^*x\bigr) \mathcal J_\mu\bigl(y^*y\bigr) \,=\, 
\frac{1}{\kappa_\mu}\int_{D_q} \mathcal J_\mu\bigl(x^*x + y^*y + x^*uvw^*y  + y^*wv^*u^*x
\bigr)\Delta(I-v^*v)^{\mu- \varrho}dv.\]
Under the coordinate transform $v\mapsto  uvw^* =:\widetilde v$ we have
$\,dv = d\widetilde v, \,\Delta(I-v^*v) = \Delta(I-\widetilde v^*
\widetilde v)$ and $\,v\in D_q\Leftrightarrow\widetilde v\in D_q$. This implies \eqref{extended}.
\end{remark}

\begin{proof}[Proof of Theorem \ref{main1}]
Let $W:= \{\mu\in \b C: \mathfrak R\mu > \varrho-1\}$. By
 the asymptotic properties of the usual gamma function we have
\begin{equation}\label{kappaestim}
\kappa_\mu \sim \,\Bigl(\frac{\pi}{\mu}\bigr)^{dq^2/2} \quad \text{uniformly in $W$ as $\mu\to\infty$.}
\end{equation}
Consider now the claimed product formula \eqref{prod3}.
Its left-hand side is holomorphic and, according 
to Corollary \ref{B_estimate}, also uniformly bounded in $W$
as a function of $\mu$. In order to estimate the right-hand side
of \eqref{prod3}, note that the argument $r^2 + s^2 +rvs + sv^*r$ is positive
semidefinite for all $v\in D_q$. Moreover,
\begin{align}\label{Testim} \text{tr}(r^2 + s^2 +rvs + sv^*r) \,=\,& 
\|r\|^2 + \|s\|^2 + 2(r\vert vs)\notag\\
\leq\,&\|r\|^2 + 2\|r\|\|vs\| + \|s\|^2\,
\leq\,(\|r\|+ \|s\|)^2
\end{align}
where again $\|\,.\,\|$ is the Frobenius norm on $M_q$. For the last estimate, it was used that $I-v^*v$ is positive-definite and therefore 
\[ \|s\|^2 - \|vs\|^2 = \,\text{tr}\bigl(s(I-v^*v)s\bigr) \geq 0.\]
Thus by Corollary \ref{B_estimate},
\[ \big|\mathcal J_\mu\bigl(r^2 + s^2 +rvs + sv^*r\bigr)\big| \,\leq\,
    e^{(\|r\|+ \|s\|)^2} \quad \forall v\in D_q.\]
This easily implies that the right-hand side in \eqref{prod3} is also holomorphic
as a function of $\mu$ in $W$, and can be estimated according to
\begin{align}\Big| \frac{1}{\kappa_\mu}\int_{D_q} &
\mathcal J_\mu\bigl(r^2 + s^2 +rvs + sv^*r\bigr)\,\Delta(I-v^*v)^{\mu- \varrho}dv \Big| \notag\\
\leq\,& C\cdot\frac{1}{|\kappa_\mu|} \cdot
\int_{D_q} \Delta(I-v^*v)^{\mathfrak R \mu -\varrho}dv \,=\,
 C\cdot \frac{\kappa_{\mathfrak R \mu}}{|\kappa_\mu|}\end{align}
with a constant $C>0$ independent of $\mu$. In view of \eqref{kappaestim}, the
last expression is of the form
\[ O\bigl(|\mu|^{dq^2/2}\bigr) \quad \text{as $\mu\to\infty$ in $W$}.\]
Now define $\,f(z):=g\bigl((z+2q)\frac{d}{2}\bigr)$ where
\[ g(\mu):= \mathcal J_\mu(r^2)\mathcal J_\mu(s^2) - \frac{1}{\kappa_\mu}
\int_{D_q} \mathcal J_\mu(r^2 + s^2 +  rvs + sv^*\!r\bigr)\,\Delta(I-v^*v)^{\mu- \varrho}dv, \, \mu\in W.\]
Then $f(z) = 0$ for all $z\in \b N_0$ by Proposition \ref{convodisc}.  The above considerations further show that $f$ is holomorphic and of polynomial growth on $\mathfrak R z >-1$.
With Carlson's theorem, validity of  \eqref{prod3} follows as claimed.
\end{proof}

\begin{remarks}\label{Bochner}
1. We may as well establish a 
Bochner type integral representation for the Bessel function $J_\mathcal\mu$ with $\mathfrak R \mu > \varrho -1$ by analytic continuation.
 Indeed, for $\mu = \frac{pd}{2}$ with an integer $p\geq q$, we obtain from \eqref{intrepfaraut} that for all $x\in M_q$, 
\[ \mathcal J_\mu\bigl(x^*x\bigr)\,=\, \int_{\Sigma_{p,q}} e^{-2i(\sigma\vert \sigma_0 x)}d\sigma \,=\, \int_{\Sigma_{p,q}} e^{-2i(\widetilde\sigma\vert x)}d\sigma.\]
 If $p\geq 2q$ then according to Corollary \ref{red}  this can be written as 
\begin{equation}\label{Bochnerneu}
\mathcal J_\mu(x^*x) = \,\frac{1}{\kappa_\mu} \int_{D_q} e^{-2i(v\vert x)} \Delta(I-v^*v)^{\mu-\varrho}dv.\end{equation}

Analytic continuation with respect to $\mu$ as above shows that 
\eqref{Bochnerneu} remains valid for all $\mu\in \b C$ with $\mathfrak R \mu >\varrho-1$. 
From this identity, it follows by the Riemann-Lebesgue Lemma for the additive matrix group $(M_q, +)$  that $\mathcal J_\mu\in C_0(\Pi_q)$. 
In case $\b F=\b R$, the 
integral representation \eqref{Bochnerneu} for $\mathfrak R \mu >\varrho-1$ goes back to \cite{H}, where it was proven  by a different method.
\medskip

\noindent
2. In the rank one case $q=1$ we have
\[ \mathcal J_\mu\bigl(\frac{r^2}{4}\bigr) = j_{\mu-1}(r) \quad (r\in \b R_+)\]
with the one variable Bessel functions
\[ j_\alpha(z)\,=\,_0F_1\bigl(\alpha +1; -\frac{z^2}{4}\bigr) = 
\sum_{n=0}^\infty
\frac{(-1)^n \Gamma(\alpha+1)}{n!\,\Gamma(\alpha+n+1)}\bigl(\frac{z}{2}\bigr)^{2n}\,.\]
Formulas \eqref{Bochnerneu} and \eqref{prod3} reduce to the well-known
Mehler formula 
\[ j_\alpha(r) = c_\alpha \int_{-1}^1 e^{-irt}(1-t^2)^{\alpha-1/2}dt\]
and the product formula 
\[ j_\alpha(r)j_\alpha(s) = c_\alpha\int_{-1}^1 j_\alpha\bigl(\sqrt{r^2 + s^2 + 2rst}\,\bigr) (1-t^2)^{\alpha-1/2}dt,\]
which are both valid for $\alpha\in \b C$ with  $\mathfrak R \alpha >-1/2$.
\end{remarks}

We finish this section with an alternative form of
product formula \eqref{prod3} and the integral representation \eqref{Bochnerneu}
for the Bessel functions $\mathcal J_{\mu}$. Let
\[ B:= \{z\in M_{1,q}(\b F): |z|<1\}\]
denote the unit ball in $M_{1,q}(\b F)\cong \b F^q$ with respect to the standard norm $\,|z| = \bigl(\sum_{j=1}^q \overline{z_j}z_j\bigr)^{1/2}$. It is easily checked that
\[\Delta(I-z^*z) = 1-|z|^2.\]

\begin{lemma}\label{polarred}
The mapping
\[ P(y_1, \ldots, y_q):= \begin{pmatrix}y_1\\y_2(I-y_1^*y_1)^{1/2}\\ 
\vdots\\
  y_q(I-y_{q-1}^*y_{q-1})^{1/2}\cdots (I-y_{1}^*y_{1})^{1/2}\end{pmatrix}; \, \,\,y_1,\ldots,y_q \in B\]
establishes a diffeomorphism from $B^q$ onto $D_q$ with the Jacobi determinant
\[ |\det  dP(y)|\,=\,\prod_{j=1}^{q-1} (1-|y_j|^2)^{d(q-j)/2}.\]
\end{lemma}

\begin{proof}
For $k$ with $2\leq k\leq q$,  the mapping 
\[\, B\times D_{k-1} \longrightarrow\, D_{k}, \quad (z,w)\mapsto \, \begin{pmatrix}z\\ w\sqrt{I-z^*z}\end{pmatrix}\]
is obviously a diffeomorphism. By Lemma 2 of \cite{FT} and Fubini's theorem, its Jacobi determinant is given by $\,(1-|z|^2)^{d(k-1)/2}.$ 
$(q-1)$-fold iteration of this decomposition yields the assertion.
\end{proof}

\noindent
It is easily checked that for $v=P(y)\in D_q$,
\[ \Delta(I-v^*v) \,=\, \prod_{j=1}^q \Delta(I-y_j^*y_j) \,=\, \prod_{j=1}^q (1-|y_j|^2).\]
This implies

\begin{corollary}\label{integralD}
Let $\zeta\in \b C$ with $\mathfrak R\zeta >-1$. Then the image measure
of $\, \Delta(I-v^*v)^\zeta dv$ under $P^{-1}: D_q\to B^q$ is given by
$\,\prod_{j=1}^{q}(1-|y_j|^2)^{\zeta+d(q-j)/2}dy.$
\end{corollary}

Application to the product formula \eqref{prod3} and the integral representation \eqref{Bochnerneu} for the Bessel function
$\mathcal J_\mu$ gives

\begin{corollary}\label{alternat} Let $\mu\in \b C$ with $\mathfrak R\mu >\varrho-1$. Then 
for all $r,s\in \Pi_q$, 
\begin{align}\mathcal J_\mu(r^2)\mathcal J_\mu(s^2) \,=\,& \frac{1}{\kappa_\mu}\int_{B^q}
\mathcal J_\mu\bigl(r^2 + s^2 +rP(y)s + sP(y)^*r\bigr) \prod_{j=1}^{q} (1-|y_j|^2)^{\mu-\varrho+d(q-j)/2}dy;\notag\\
\mathcal J_\mu(r^2) = \,&\frac{1}{\kappa_\mu}\int_{B^q}e^{-2i(P(y)\vert r)}
\prod_{j=1}^q (1-|y_j|^2)^{\mu-\varrho+d(q-j)/2}dy.\notag
\end{align}
\end{corollary}

\subsection{Bessel convolutions on the cone $\Pi_q$}\label{cone_hypergroup}

For real $\mu >\varrho-1$ 
 the measure $\kappa_\mu^{-1}\cdot\Delta(I-v^*v)^{\mu-\varrho}$ in  product formula \eqref{prod3}
is a probability measure on $D_q$. We shall see that \eqref{prod3} leads to a hypergroup
convolution on $\Pi_q$ which is of the same form as those of Proposition \ref{convodisc}.
This will give us a continuous series of commutative hypergroup structures
on $\Pi_q$ which interpolate those occuring as orbit hypergroups for the indices $\mu= pd/2$, $p\geq 2q$ an integer.

\begin{theorem}\label{main2}
Fix an index $\mu\in\b R$ with $\mu > \varrho-1.$ 
\begin{enumerate}
\item[\rm{(a)}] The assignment
\[(\delta_r *_\mu \delta_s)(f) := \frac{1}{\kappa_\mu}\int_{D_q} f\bigl(\sqrt{r^2 + s^2 + rvs + sv^*r}\,\bigr)\,
\Delta(I-v^*v)^{\mu-\varrho}\, dv; \quad f\in C(\Pi_q)\]
defines a commutative hypergroup structure $\Pi_{q,\mu} = (\Pi_q\,, *_\mu)$ with neutral element $0$
and the identity mapping as involution. The support of $\delta_r*_\mu\delta_s$ satisfies 
\[\text{supp}(\delta_r*_\mu\delta_s) \subseteq  \{t\in \Pi_q: \|t\|\leq \|r\|+\|s\|\}.\]
\item[\rm{(b)}] A Haar measure of the hypergroup $\Pi_{q,\mu}$ is given by
\[ \omega_\mu(f) = \frac{2^{-q\mu}}{\Gamma_{\Omega_q}(\mu)}\int_{\Omega_q} f(\sqrt{r}) \Delta(r)^\gamma dr\]
with $\,\gamma = \mu-\frac{d}{2}(q-1)-1\,=\, \mu-\frac{n}{q}.$
\end{enumerate}
\end{theorem}

\begin{remark} The specific normalization of the Haar measure is motivated by Theorem \ref{dualspace} below.
\end{remark}

\begin{proof}
Ad (a). Clearly  $\delta_r *_\mu \delta_s$ is a probability measure on $\Pi_q$ and 
\[\text{supp}(\delta_r*_\mu\delta_s) = \big\{ \sqrt{r^2 + s^2 + rvs + sv^* v}\,: v\in \overline{D_q}\big\},\]
which does not depend on $\mu$. 
The stated support inclusion is immediate from estimate 
\eqref{Testim}. This shows property (2) in the definition of a hypergroup. (3) is clear because it is known to be true for indices $\mu= pd/2$ which 
lead to an orbit hypergroup structure. Property (4) is obvious, and (5) with $\overline{r} = r$ is again true in general because it is true in 
 the orbit hypergroup cases. For the proof of (1) it suffices
to show that for each $f\in C_b(\Pi_q)$ the mapping $(r,s)\mapsto f(r*_\mu s)$ is continuous. But this is clear from the continuity of the map 
$\, (r,s,v)\mapsto f\bigl(\sqrt{r^2+s^2 + rvs + sv^*\!r}\bigr)$ on $\Pi_q^{\,2}\times D_q$. Commutativity of the convolution is clear, and for the 
proof of associativity it again suffices to consider point measures.
So let $r,s,t\in \Pi_q$ and $f\in C_b(\Pi_q)$. Then
\[\delta_r*_\mu(\delta_s*_\mu\delta_t)(f) = 
\frac{1}{\kappa_\mu^2}\int_{D_q}\int_{D_q} f\bigl(H(r,s,t;v,w)\bigr)
\Delta(I-v^*v)^{\mu-\varrho}\Delta(I-w^*w)^{\mu-\varrho}dv dw \,=:I(\mu) \]
with a certain $\Pi_q$-valued argument $H$ that is
independent of the index $\mu$. Similar, 
\[(\delta_r*_\mu\delta_s)*_\mu\delta_t(f) = I^\prime(\mu)\]
with some $\mu$-independent argument $H^\prime$ instead of $H$. 
The integrals $I(\mu)$ and $I^\prime(\mu)$ are well-defined and holomorphic in 
$\{\mu\in \b C: \mathfrak R\mu >\varrho -1\}$. 
Further, we know that $\,I(\mu) = I^\prime(\mu)\,$ 
 for all $\mu = pd/2$ with an integer $p\geq 2q$. By analytic continuation as in the proof of Theorem \ref{main1} (use \eqref{kappaestim} again)
we obtain validity of this relation for all $\mu$ with $\mathfrak R \mu >\varrho-1$.  

\noindent
Ad (b). We have to prove that
\[ \int_{\Omega_q} f(s*_\mu \sqrt{r}) \Delta(r)^\gamma dr\,=\, \int_{\Omega_q} f(\sqrt{r})
\Delta(r)^\gamma dr \]
for all $f\in C_c(\Pi_q)$ and $s\in \Pi_q$. Equivalently,
\begin{align}\label{haarrel}
\int_{\Omega_q}\int_{D_q} f\Bigl(&\sqrt{s^2 + r + sv\sqrt{r}\, + \sqrt{r}v^*s}\,\Bigr)\Delta(I-v^*v)^{\mu-\varrho} \Delta(r)^{\mu-\frac{n}{q}} dr dv \notag\\ 
&\,=\,\int_{\Omega_q} f(\sqrt{r})\Delta(r)^{\mu-\frac{n}{q}}dr.
\end{align}
We know from Section \ref{orbitcone} that this relation is true if $\mu=pd/2$ with an integer $p\geq 2q$. Moreover, the stated result is well known in case $d=q=1$
for all $\mu$ in question; actually, in these cases
$*_\mu$ is just the convolution of a Bessel-Kingman hypergroup (see Section 3.5 of \cite{BH}). So we may assume that $dq\geq 2$.
Let again $\,W= \{\mu\in \b C: \mathfrak R\mu > \varrho -1\}$. Then $\mathfrak R\mu - \frac{n}{q} > 0$ for all $\mu\in W$, and for fixed $f$ both sides of \eqref{haarrel} are 
well-defined and holomorphic as functions of $\mu$ in $W$.
We shall again carry out analytic continuation with respect to $\mu$, based on Carlson's theorem. We proceed in two steps.

Step 1. For $R>0$ let $\, B_R:= \{ r\in \Pi_q: \|r\| \leq R\}$. 
Choose $R>0$ with $R^q \leq e^{\pi/d}$. The explicit formula for the convolution 
$*_\mu$ shows that there exist constants $S_0, R_0 >0$ such that
for $f\in C_c(\Pi_q)$ with $\,\text{supp}\, f\subseteq B_{R_0}$ and for $s\in \Pi_q$ with $\|s\| \leq S_0$ we have
\[ f(s*_\mu \sqrt{r}\,) = 0 \quad\text{if }\, \|r\|\geq R.\]
Suppose that $|f|\leq C$ on $\Pi_q$ and notice that $\Delta(r)  \leq \|r\|^q$. 
The left-hand side of \eqref{haarrel} may then be estimated according to

\begin{align}
|LHS| \leq \,& C\int_{B_R}\Delta(r)^{\mathfrak R \mu-\frac{n}{q}} \,dr\cdot\int_{D_q} \Delta(I-v^*v)^{\mathfrak R \mu - \varrho} \,dv \notag\\
\leq \,&C\cdot \kappa_{\mathfrak R \mu} \cdot\text{vol}(B_R)\cdot (R^q)^{\mathfrak R \mu-\frac{n}{q}}
\end{align}
Due to the asymptotics \eqref{kappaestim}, $\mu\mapsto \kappa_{\mathfrak R \mu}$ is bounded on $W$.  Thus by our  initial assumption on $R$ we arrive
at
\[ |LHS| = O\bigl(e^{\frac{\pi}{d}|\mu|}\bigr) \quad\text{as $\mu\to\infty$ in $W$}.\]
The right-hand side of 
\eqref{haarrel} coincides with the left-hand side
for $s=0$ and is therefore of (at most) the same order as a function of $\mu$.
As in the proof of Theorem \eqref{main1}, we have to substitute $\mu = (z+2q)\frac{d}{2}$
to obtain holomorphic functions in $z$ on $\mathfrak R z>-1$ which coincide for all $z\in \b N_0$ and are of order $\, O\bigl(e^{\pi|z|/2}\bigr)$ in $\mathfrak R z \geq 0$. 
Carlson's theorem now implies the assertion under the above restrictions on $f$ and $s$.

\smallskip

Step 2. Let $f\in C_c(\Pi_q)$ and $s\in \Pi_q$ be arbitrary. For $\delta >0$ define 
$f_\delta(r):= f(\delta r)$. We shall use that the convolution $*_\mu$ is homothetic, i.e. 
\[f_\delta (r*_\mu s) = f(\delta r *_\mu \delta s).\]
Fix constants $S_0, \, R_0$ as in step 1 and choose $\delta >1$ such that $\text{supp}f_\delta\subseteq B_{R_0}$ and $\big\|\frac{s}{\delta}\big\|\leq S_0$. Then
\begin{align} 
I(s):= \,&\int_{\Omega_q} f(s *_\mu \sqrt{r}) \Delta(r)^\gamma dr \,=\,\int_{\Omega_q} 
f_\delta\Bigl(\frac{s}{\delta} *_\mu \frac{\sqrt{r}}{\delta}\Bigr) \Delta(r)^\gamma dr\notag\\
=\,& \delta^{2q\gamma +2n} \int_{\Omega_q} 
f_\delta\Bigl(\frac{s}{\delta} *_\mu \sqrt{t}\Bigr)\Delta(t)^\gamma dt. \notag
\end{align}
As $f_\delta$ and $\frac{s}{\delta}$ satisfy the conditions of step 1, we 
conclude that $I(s) = I(0)$, which 
 finishes the proof of part (b).
\end{proof}

Our next aim is to determine the dual and the Plancherel measure 
 of the hypergroup $\Pi_{q,\mu}$ with $\mu > \varrho-1$. 
For $s\in \Pi_q$, we define
\[ \phi_s(r) = \phi_s^\mu(r):= \mathcal J_\mu\bigl(\frac{1}{4}sr^2s\bigr), \quad r\in \Pi_q\,.\]
Notice that 
\begin{equation}\label{symm}\phi_s(r) = \phi_r(s) \end{equation}
because $J_\mu$ depends only on the eigenvalues of its argument. Moreover, we have

\begin{lemma}\label{vanish2} Let $\mu>\varrho-1$. Then for each $s\in \Pi_q$, $\, \phi_s = \phi_s^\mu$ belongs to $C_b(\Pi_q)$ with $\|\phi_s\|_\infty = \phi_s(0) = 1$. If $s\in \Omega_q$, then even
 $\, \phi_s\in C_0(\Pi_q) $.
\end{lemma}

\begin{proof} The first assertion is immediate from the Bochner-type integral representation \eqref{Bochnerneu} for $\mathcal J_\mu$. Now suppose $s\in \Omega_q$. Then $r\to\infty $ in $\Pi_q$ implies that $r^2\to\infty $ and also $sr^2s\to\infty $, because  $s$ is invertible. The second assertion thus follows from 
the fact that $\mathcal J_\mu$ vanishes at infinity.
\end{proof}

From this Lemma together with the product formula \eqref{extended}  it is immediate
that each $\phi_s$ with $s\in \Pi_q$ belongs to $\widehat\Pi_{q,\mu}.$
The results of Section \ref{orbitcone} suggest that these  Bessel functions  actually make up the complete dual: 

\begin{theorem}\label{dualspace} (1) The dual space of $\Pi_{q,\mu}$ with $\mu >\rho -1$ is given by
\[\widehat\Pi_{q,\mu} = \,\{\phi_s = \phi_s^\mu : s\in \Pi_q\}.\]
(2) The hypergroup $\Pi_{q,\mu}$ is self-dual via the homeomorphism $\,\Psi:\Pi_{q,\mu} \to \widehat\Pi_{q,\mu}\,, \, s\mapsto \phi_s$.
Under this identification,  the Plancherel measure $\pi_\mu$ on $\widehat\Pi_{q,\mu}$ coincides with the Haar measure $\omega_\mu$.
\end{theorem}

For the proof of part (1) we  need the following

\begin{lemma}\label{closed}
The set $A:= \{\phi_s: \,s\in \Pi_q\}$ is closed in $\widehat\Pi_{q,\mu}$ with respect to the
topology of locally uniform convergence.
\end{lemma}

\begin{proof}
Let $(\phi_{s_j})_{j\in \b N}$ be a sequence in $A$ converging to $\alpha\in \widehat\Pi_{q,\mu}$ locally uniformly. If the sequence $(s_j) \subset \Pi_q$ is bounded, then after passing to
a subsequence  we may assume that $s_j \to s\in \Pi_q$ as $j\to\infty$.
Then $\phi_{s_j} \to \phi_s$ and therefore $\alpha = \phi_s\in A.$ If else the original sequence $(s_j)$ is unbounded, then after passing to a subsequence we may assume that $s_j\to\infty$. Thus by Lemma \ref{vanish2}, 
$\phi_{s_j}(r) = \phi_r(s_j) \to 0$ for
all $r\in \Omega_q$ as $j\to\infty$. This implies $\alpha =0$ which contradicts the
convention
$0\notin\widehat\Pi_{q,\mu}$. 
\end{proof}

\begin{proof}[Proof of Theorem \ref{dualspace}] In a first step, we establish that  $\Pi_{q,\mu}$ has subexponential growth
in the sense of \cite{V1}, i.e. for each compact subset $K\subset \Pi_q$ 
and each $c>1$, the Haar measure satisfies $\omega_\mu(K^{m}) = o(c^m)$. Here $K^{m}$ denotes
the $m$-fold convolution power of $K$,  the  convolution product of subsets $A,B$ of a hypergroup $(X,*)$ being defined by $A*B = \bigcup_{x\in A, y\in B} \text{supp}(\delta_x*\delta_y)$. Once subexponential growth   is known, Theorem 2.17 of \cite{V1} will imply that the support of
the Plancherel measure $\pi_\mu$  coincides with the complete dual $\widehat\Pi_{q,\mu}$. 
For the proof of subexponential growth, it suffices to consider  the balls $B_R = \{r\in \Pi_q: \|r\|\leq R\}$. From the support properties of $*_\mu$ we see that $B_R^{\,m} \subseteq B_{mR}$.
Moreover,
\[ \omega_\mu(B_R) =\, C\cdot\int_{\|\sqrt{r}\|\leq R} \Delta(r)^\gamma dr = \, O\bigl(R^{2q\gamma + 2n}\bigr) \,=\, O\bigl(R^{2q\mu}\bigr) \quad \text{as }\, R\to\infty.\]
Thus for fixed $R>0$, we obtain $\,\omega_\mu(B_R^{\,m}) =  O\bigl(m^{2q\mu}\bigr)$, and the assertion follows.

In a second step, we determine $\pi_\mu$. The decisive ingredient will be 
known results about the  Hankel transform on a symmetric cone. 
 For $s,r\in \Omega_q$, define
\[ H_\mu(s,r) = \frac{1}{\Gamma_{\Omega_q}(\mu)}\mathcal J_\mu\bigl(\sqrt{s}r\sqrt{s}\bigr).\]
Suppose that  $\mu > d(q-1) +1$. Then according to Theorem XV.4.1. of \cite{FK}, the Hankel transform
\[ U_\mu F(s) := \int_{\Omega_q} H_\mu(s,r)F(r) \Delta(r)^\gamma dr\]
defines an isometric and involutive isomorphism of $L^2(\Omega_q, \Delta(r)^\gamma dr)$. The argumentation of \cite{H} and \cite{FT}, Section 5 shows that this statement actually extends to all $\mu\in \b R$ with $\mu > \frac{d}{2}(q-1)$, i.e.
$\gamma >-1$. Let $f\in L^2(\Pi_q, \omega_\mu)$. Then $F(r):= f\bigl(\sqrt{r}\bigr)$ belongs to $L^2(\Omega_q, \Delta(r)^\gamma dr)$ and a short calculation shows that
\[\widehat f(\phi_s) = \int_{\Pi_q} \phi_s(r) f(r) d\omega_\mu(r)\,=\,
2^{-q\mu}\, U_\mu F\bigl(\frac{s^2}{4}\bigr) \quad\forall\,s\in \Omega_q.\]
Moreover, by the isometry of $U_\mu$ we readily obtain
\[ \int_{\Pi_q}|\widehat f(\phi_s)|^2 d\omega_\mu(s) \,=\, \int_{\Pi_q}|f(s)|^2d\omega_\mu(s).\]
This shows that the Plancherel measure $\pi_\mu$ associated with $\omega_\mu$ is given by
\[ \pi_\mu(g) = \int_{\Pi_q} g(\phi_s) d\omega_\mu(s), \quad g\in C_c(\widehat\Pi_{q,\mu}).\]
Hence the support of $\pi_\mu$, which we already know to coincide with $\widehat \Pi_{q,\mu}$, also coincides with the closure of 
the set $\{\phi_s: \,s\in \Pi_q\}$ in $\widehat\Pi_{q,\mu}$ with respect to the
topology of locally uniform convergence. The proof of part (1) is therefore
accomplished by Lemma \ref{closed}. For part (2), it remains to verify that $\Psi$ is a homeomorphism. Continuity and surjectivity are clear. For injectivity, 
suppose $\phi_s = \phi_r$. Then in view of \eqref{symm} we have $\widehat\delta_s = \widehat \delta_r$ and the injectivity of the Fourier transform of measures on the hypergroup  $\Pi_{q,\mu}$ implies $s=r$. 
To check continuity of $\Psi^{-1}$ suppose that $\phi_{s_i} \to \phi_s$ locally uniformly. Then $\, \widehat\delta_{s_i} \to 
\widehat \delta_{s}$ locally uniformly on $\widehat\Pi_{q,\mu}$. Levy's continuity theorem (Thm. 4.2.2. in \cite{BH}) imlies that $\delta_{s_i} \to \delta_s$ weakly, and hence $s_i\to s$.
\end{proof}

\subsection{The limit case $\mu=\varrho-1$.}

Using Corollary \ref{integralD} we see that 
the convolution $*_\mu$ with $\mu > \varrho -1$ can be written 
in the alternative form 
\begin{equation} \label{alt2}(\delta_r*_\mu\delta_s)(f) = \frac{1}{\kappa_\mu} 
\int_{B^q} f\bigl(\sqrt{r^2 + s^2 + rP(y)s + sP(y)^*r}\,\bigr) \prod_{j=1}^q(1-|y_j|^2)^{\mu-\varrho +\frac{d}{2}(q-j)}dy.\end{equation}
We shall use this representation  to determine the limit of the  convolution $*_\mu$
as $\mu\downarrow\rho -1$,  where it assumes a degenerate form. As $\rho -1 = pd/2$ with 
$p = 2q-1$, it is natural to expect that the reulting limit convolution  coincides
with the orbit hypergroup convolution $*_{\rho -1} $ on $\Pi_q$ derived from $M_{q,2q-1}.$ 
In the following, $d\sigma$ denotes the normalized surface measure on the unit
sphere $S= \{z\in M_{1,q}: |z|=1\}$. 

The coordinate transform $P:\,B^q\to D_q$ of Lemma \ref{polarred} is assumed to be continuously extended to  $\overline B^{\,q}$. 

\begin{proposition}\label{limitcase}
As $\mu\downarrow \varrho-1$, the convolution product $\delta_r*_\mu\delta_s$ converges weakly to the probability measure $\delta_r\,\tilde *\,\delta_s$ on $\Pi_q$  given by
\begin{align} (\delta_r \,\tilde *\,\delta_s)(f) = \,\widetilde\kappa
\int_{B^{q-1}}\!\int_S
 f\bigl(\sqrt{r^2 + s^2 + rP(y)s + sP(y)^*r}\,\bigr) \prod_{j=1}^{q-1}
(1-&|y_j|^2)^{\mu-\varrho +\frac{d}{2}(q-j)}\cdot\notag\\
&dy_1\ldots dy_{q-1}\, d\sigma(y_q)\notag
\end{align}
with a normalization constant $\,\widetilde\kappa > 0$. 
The product $\,\tilde *\,$ defines a commutative hypergroup structure 
on $\Pi_q$ which coincides with the orbit hypergroup $\Pi_q^{\rho -1}$ derived from $M_{q,2q-1}$ as in Section 
\ref{orbitcone}. In particular, 
\[ (\delta_r \,\tilde *\,\delta_s)(f)\,= \,(\delta_r *_{\rho -1}\!\delta_s)(f)
= \,\int_{\Sigma_{q,2q-1}} f\bigl(\sqrt{r^2 + s^2 + r\widetilde\sigma s +(r\widetilde\sigma s)^*}\bigr) d\sigma\]
and the additional statements of Theorem \ref{main2}
and Theorem \ref{dualspace} extend to the case
$\mu = \varrho-1$. 
\end{proposition}

\begin{proof}  
For $\mu>\varrho-1$ consider the probability measure 
\[ p_\mu:= c_\mu^{-1} (1-|y|^2)^{\mu-\varrho}\,1_{B}(y)dy\]
on $M_{1,q}$, where $c_\mu:= \int_B (1-|y|^2)^{\mu-\varrho} dy \,$ 
and  $1_B$ denotes the characteristic function of the ball $B$. It is
easily checked that  $p_\mu$ tends weakly to the normalized
surface measure $d\sigma$ on $S$ as $\mu\to\varrho-1$. Indeed, let $f\in C(M_{1,q})$ and put
$F(\tau) := \int_S f(\tau y)d\sigma(y), \, \tau \geq 0.$ Then with $c_\mu^\prime = \int_0^1 (1-\tau^2)^{\mu-\varrho} \tau^{dq-1}d\tau$,
\[ \int_B f dp_\mu\,=\,\frac{1}{c_\mu^\prime} \int_0^1 F(\tau)(1-\tau^2)^{\mu-\varrho} \tau^{dq-1}d\tau\,\longrightarrow F(1) \quad\text{as }\, \mu\to\varrho-1.\]
This proves that $\delta_r*_\mu\delta_s \to \delta_r\tilde *\,\delta_s$ 
weakly. It is clear that $\, \text{supp}(\delta_r\tilde *\,\delta_s)\subseteq
\,\text{supp}(\delta_r*_\mu\!\delta_s)$ for $\mu >\varrho-1$ which implies
the same support inclusion as in Theorem \ref{main2} (a). In the limit $\mu \to \varrho -1$ we further obtain that the Bessel functions 
\[\phi_s(r) = \mathcal J_{\varrho-1}\bigl( \frac{1}{4}sr^2s\bigr), \quad s,r\in \Pi_q\]
satisfy the product formula
\[ \phi_s(r) \phi_s(t) \, = \, \int_{\Pi_q}
\phi_s(\tau) d(\delta_r\tilde *\,\delta_t)(\tau) \quad \forall\, s\in \Pi_q.\] 
On the other hand, 
consider the orbit hypergroup  $\,X_{\varrho -1}=(\Pi_q, *_{\varrho -1})$ 
derived from $M_{p,q}$ with $p=2q-1$.  
Its dual space consists exactly of the Bessel functions $\phi_s(r) = \mathcal J_{\rho-1}(\frac{1}{4}rs^2r)$ as above. In particular, 
\[ \phi_s(r) \phi_s(t) \, = \, \int_{\Pi_q}
\phi_s(\tau) d(\delta_r*_{\!\varrho-1}\delta_t)(\tau) \quad \forall\, s\in \Pi_q.\] 
The injectivity of the Fourier transform on the hypergroup $X_{\varrho -1}$ now implies
that $\delta_r*_{\!\varrho-1}\delta_t\,=\, \delta_r\tilde *\,\delta_t$ for all $r,t$.\end{proof}

\begin{remarks}
1. We conjecture that the formulas of Corollary \ref{alternat} permit degenerate
extensions to successively larger index ranges. As soon as the exponent 
in one of the iterated integrals becomes critical, the corresponding integral 
over  $B$ should be replaced by an integral over $S$. More precisely, we conjecture that within the range
\[\big\{\mu\in \b R: \varrho -\frac{d}{2}(q-k) -1   < \mu    \leq \varrho -\frac{d}{2}(q-k-1) -1\big\}, \quad k= q-1,\ldots, 1\]
the following product formula is valid: 
\begin{align}\mathcal J_\mu(r^2) \mathcal J_\mu(s^2) \,=\,\frac{1}{\kappa_{\mu,j}} \int_{B^k} \int_{S^{q-k}}
 \mathcal J_\mu\bigl(r^2 + s^2 + rP(y)s +& sP(y)^*r\bigr) 
\prod_{j=1}^k(1-|y_j|^2)^{\mu-\varrho + \frac{d}{2}(q-j)}\notag\\& \cdot dy_1\ldots dy_k\, d\sigma(y_{k+1})\dots d\sigma(y_q).\notag
\end{align}
Also, there should be analogous integral representations for $\mathcal J_\mu$
within the above ranges of $\mu$.

\medskip
2. Further properties of the hypergroups $\Pi_{q,\mu}$, concerning 
their automorphism groups as well as stochastic aspects (such as limit theorems for random walks on matrix cones associated with Bessel convolutions)
will be the subject of a forthcoming paper joint with M. Voit. 
\end{remarks}

\section{Hypergroups associated with rational Dunkl operators of type $B$}

In the analysis of the previous sections, one may be interested in
questions which depend only on the spectra of the matrices from the underlying cone
$\Pi_q = \Pi_q(\b F)$. This amounts to considering functions and measures on $\Pi_q$ which are invariant under unitary
conjugation. 
For $x\in H_q$ we denote by $\sigma(x) = (\xi_1,\ldots , \xi_q)\in \b R^q$ the set of eigenvalues of $x$ ordered by size, i.e. $\xi_1\geq \ldots \geq \xi_q$. 
The unitary group $U_q$ acts on $\Pi_q$ via conjugation, $(u,r)\mapsto uru^{-1}$, and 
the orbits under this action
are parametrized
by the set $\Xi_q$ of possible spectra $\sigma(r)$
of matrices $r\in \Pi_q$, 
\[\Xi_q = \{\xi = (\xi_1, \ldots \xi_q)\in \b R^q: 
\xi_1\geq\ldots\geq \xi_q\geq 0\}.\]
We mention that $U_q$ is the maximal compact subgroup of the automorphism group
of $\Omega_q$; inded the latter coincides with $GL_q = GL_q(\b F)$, acting on $H_q$ via $(g,x)\mapsto gxg^*$. The set $\Xi_q$ is a closed Weyl chamber of the hyperoctahedral group $B_q = 
S_q\ltimes \b Z_2^q$ which acts on $\b R^q$ by permutations of the basis vectors and sign changes.
In Section \ref{convocones} we saw that $\Pi_q$ carries a continuously parametrized family of commutative hypergroup 
structures $*_\mu$ with $\mu\geq\varrho-1$, as well as additional orbit hypergroup structures for $\mu = pd/2, \, p\geq q$ an integer.   In the following section we are going  to show that under the above action of $U_q$ on $\Pi_q$ each convolution $*_\mu$  induces 
a commutative  orbit hypergroup convolution $\circ_\mu$ on $\Xi_q$, similar as orbit hypergroups were obtained from (commutative) groups in Section \ref{orbithypergroups}.
In Section \ref{dunklbessel} we shall then identify the characters of $(\Xi_q, \circ_\mu)$
with multivariable Bessel functions of Dunkl type which are
associated with the root system of type $B_q$. In effect, we thus obtain 
a continuous series of commutative hypergroup structures on the chamber $\Xi_q$ whose characters
are given by Dunkl-type Bessel functions.

\subsection{Convolutions on the spectra of positive definite matrices}\label{spectra}\label{chamberconvo}

In the situation described above, the canonical mapping
\[ \pi: \Pi_q \to \Xi_q\,, \quad r\mapsto \sigma(r)\]
is continuous, surjective and open with respect to the standard topologies on both sets.
Therefore the map $\, \Pi_q^{U_q} \to \Xi_q, \, U_q .r\mapsto \sigma(r)\,$
becomes a homeomorphism when $\Pi_q^{U_q}$ is equipped with the quotient topology. 
In the following, $du$ denotes the normalized Haar measure on $U_q$ and $\xi\in \Xi_q$ is always identified with the diagonal matrix $\text{diag}(\xi_1,\ldots, \xi_q)\in \Pi_q$ without mentioning. Moreover, we introduce the index set
\[ \mathcal M_q:= \big\{\frac{pd}{2}, \, p = q, q+1, \ldots\big\} \cup \, ]\rho -1,\infty[\,.\]

\begin{theorem}\label{main3} (1) For each $\mu \in \mathcal M_q$
the chamber $\Xi_q$ carries a commutative hypergroup structure with convolution
\[
 (\delta_\xi \circ_\mu \delta_\eta)(f) := \,\int_{U_q} (f\circ\pi)(\xi *_\mu u\eta u^{-1}) du, \quad f\in C(\Xi_q).\]
The neutral element of the hypergroup $\,\Xi_{q,\mu}:=(\Xi_q,\circ_\mu)$ is  $\,0\in \Xi_q$ and the involution is given by the identity mapping.

\noindent (2) A Haar measure on $\Xi_{q,\mu}$ is given by
\[
\widetilde\omega_\mu \,= \,\pi(\omega_\mu) \,=\,d_\mu h_\mu(\xi)d\xi \quad\text{with}\quad
 h_\mu(\xi) =\,\prod_{i=1}^q \xi_i^{2\gamma + 1} \prod_{i<j} (\xi_i^2-\xi_j^2)^d\]
and a constant  $ d_\mu >0.$
\end{theorem}

\begin{remarks}\label{remmain3}
1. The constant $d_\mu$ will be determined in Section \ref{dunklbessel}.

\noindent
2. In the generic case $\mu>\varrho-1$ the convolution $\circ_\mu$ can be more explicitly written as
\[(\delta_\xi \circ_\mu \delta_\eta)(f)\,=\, \frac{1}{\kappa_\mu}
\int_{D_q} \int_{U_q} f\bigl(\sigma\bigl(\sqrt{\xi^2 + u\eta^2u^{-1} + \xi v u\eta u^{-1} + u\eta u^{-1}
v^*\xi}\, \bigr)\bigr)\Delta(I-v^*v)^{\mu-\varrho}dudv.\]
\end{remarks}

\noindent
We start with some preparations for the proof of Theorem \ref{main3}.

\begin{lemma}\label{l1}
For $u\in U_q$ consider the automorphism of $\Pi_q$ given by
$\, T_u: r\mapsto uru^{-1}$. 
The image measure of the convolution product $\delta_r*_\mu\delta_s\in M^1(\Pi_q)$ under $T_u$ is given by
\[ T_u(\delta_r*_\mu\delta_s) = \delta_{uru^{-1}} *_\mu \delta_{usu^{-1}}.\]
\end{lemma}

\begin{proof}
Recall that the measure $\Delta(I-v^*v)^{\mu-\varrho}dv$ in the convolution formula of Theorem \ref{main2} is 
invariant under unitary conjugations. This yields immediately that 
\[ (g\circ T_u)(r*_\mu s) \,=\,g(uru^{-1}\,*_\mu \,usu^{-1}) \quad \forall \,g\in C(\Pi_q). \]
\end{proof}

\noindent
Aside: this lemma just says that the automorphism $T_u$ of $\Pi_q$ 
is actually a hypergroup automorphism for each of the convolutions $*_\mu$. The next lemma is a special case of \cite{FK}, Theorem VI.2.3. 

\begin{lemma}\label{l2}
For integrable  functions $g:\Pi_q\to\b C$, 
\[ \int_{\Pi_q} g(r)dr \,=\, \kappa\!\int_{\Xi_q}\int_{U_q} g(u\xi u^{-1})\, du\prod_{i<j}(\xi_i - \xi_j)^d\, d\xi\]
with a normalization constant $\kappa=\kappa_q >0$.
\end{lemma}

\begin{proof}[Proof of Theorem \ref{main3}]
Ad (1). We shall employ the technique of \cite{Je}, Section 13  for the transfer of hypergroup 
structures via orbital mappings. Notice first that the continuous open surjection $\pi:\Pi_q\to \Xi_q$ is also proper (because $U_q$ is compact) and thus provides
an orbital mapping from the hypergroup $\Pi_{q,\mu}$ onto $\Xi_q$ in the sense of \cite{Je}, Section 13.
For $\xi\in \Xi_q$ define
\[ \epsilon_\xi :=\int_{U_q} \delta_{u\xi u^{-1}}du \]
which is a probability measure on $\Pi_q$ and satisfies $\,\text{supp}\,\epsilon_\xi = \pi^{-1}(\xi).$ We claim that each $\epsilon_\xi$ is $\pi$-consistent in the following sense:
\begin{equation}\label{consistent} \pi(\epsilon_\xi*_\mu\delta_s) = \pi(\epsilon_\xi*_\mu\delta_t) \quad\text{for all $s,t\in \Pi_q$ with $\pi(s)= \pi(t)$};\end{equation}
here $\pi$ is extended to $M^+(\Pi_q)$ by taking image measures. 
For the proof of \eqref{consistent}, suppose that $s,t\in \Pi_q$ satisfy 
$\pi(s) = \pi(t)$. Then for  each $g\in C(\Pi_q)$ which is invariant under $U_q$-conjugation we have 
\[ \int_{\Pi_q} g\,d(\epsilon_\xi *_\mu\delta_s)\,=\, \int_{U_q} g(u\xi u^{-1}*_\mu s)du \,=\, \int_{U_q} g(\xi *_\mu\, u^{-1}su)du.\]
Notice that for the second identity Lemma \ref{l1} has been used.
As $s$ and $t$ have the same spectra, the last integral does not change when $s$ is replaced by $t$. This proves \eqref{consistent}.
The orbital mapping $\pi$ also 
satisfies $\pi(0) = 0\in \Xi_q$ and $\pi^{-1}(0) =0\in \Pi_q$. We can now apply \cite{Je}, Theorem 13.5.A.  This shows that $\Xi_q$ becomes a commutative hypergroup
with convolution 
\[ \delta_\xi \circ_\mu \delta_\eta \,=\, \pi(\epsilon_\xi *_\mu \epsilon_\eta),\]
the identity mapping as involution and neutral element $\pi(0) = 0.$ 
This proves the assertions of part (1).

\smallskip\noindent
Ad (2). According to Theorem 13.3.A of \cite{Je}, 
a Haar measure $\widetilde\omega_\mu$ on $\Xi_{q,\mu}$ is given
by the image measure of $\omega_\mu$ under $\pi$. Let $f\in C_c(\Xi_q)$ and put $g = f\circ\pi$, which is $U_q$-invariant. Then
\[\int_{\Xi_q} fd\widetilde \omega_\mu\,=\,\kappa\!\int_{\Pi_q} g\,d\omega_\mu \,=\, 
\kappa^\prime \!\int_{\Xi_q} g(\sqrt{\xi}\,) \prod_{i=1}^q\xi_i^\gamma  \prod_{i<j}(\xi_i-\xi_j)^d d\xi \]
where $\,\sqrt{\xi} = (\!\sqrt{\xi_1}, \ldots, \sqrt{\xi_q}\,)$. Up to a constant factor, the last integral coincides with 
$\,\int_{\Xi_q} g(\xi) h_\mu(\xi) d\xi.$ 
\end{proof}

\begin{remarks}\label{double}
1. Recall that for $\mu=pd/2$ with an integer $p\geq q$ the hypergroup $\Pi_{q,\mu}$ is just the orbit hypergroup  obtained from the multiplication action of the unitary group $U_p$ on  $M_{p,q}.$
In this case, the above hypergroup structure of $(\Xi_q, \circ_\mu)$ can also be described
as an orbit hypergroup derived directly from $M_{p,q}$, as follows:
Consider the action of the group $L:= U_p\times U_q$ on $M_{p,q}$ by
\[ x\mapsto uxv^{-1}, \quad (u,v)\in L.\]
The orbits of this action are parametrized by the possible sets of singular values of matrices from $M_{p,q}$. Indeed, let
$\,\sigma_{sing}(x) = \sigma\bigl(\sqrt{x^*x}) = (\xi_1,\ldots , \xi_q) \in \Xi_q\,$ denote the singular spectrum of $x\in M_{p,q}$, the singular values being
ordered by size. We have the equivalences
\begin{align}\sigma_{sing}(x) = \sigma_{sing}(y) \,&\Longleftrightarrow\,
\sigma(x^*x)= \sigma(y^*y) \,\notag\\
&\Longleftrightarrow \, \exists v\in U_q: 
y^*y = vx^*xv^{-1} = (xv^{-1})^*(xv^{-1}) \notag\\
&\Longleftrightarrow \, \exists (u,v)\in U_p\times U_q:\, y=uxv^{-1}.\notag
\end{align}
Therefore the orbit space $M_{p,q}^L$ can be identified with the chamber $\Xi_q$ 
via $\, U_p\,x\,U_q \mapsto \sigma_{sing}(x)$, and this is easily checked to be
a homeomorphism with respect to the natural topologies on both spaces.
Notice that $\, (uxv^{-1})^*uxv^{-1} = vx^*xv^{-1}.$ Hence under the mapping $\,\Phi:M_{p,q}\to \Pi_q, \, x\mapsto \sqrt{x^*x},$ 
the above action of $L$ on $M_{p,q}$ induces the conjugation action of $U_q$ on 
$\Pi_q$. Moreover, the orbit convolution on $\Pi_{q,\mu}$ is defined in such a way that 
\[ \delta_r *_\mu \delta_s \,=\, \Phi(Q_r * Q_s) \quad \text{ with }\quad Q_r = \int_{U_p} \delta_{u\sigma_0r}\,du \,\in M^1(M_{p,q})\]
where $*$ denotes the usual convolution on the additive group $M_{p,q}$. This shows that for 
$\mu = pd/2$, the convolution of the hypergroup $\Xi_{q,\mu}$ coincides with the convolution of the 
orbit hypergroup 
$M_{p,q}^L \cong \Xi_q$, which is in turn naturally identified with the 
convolution of the Gelfand pair 
$(L\ltimes M_{p,q},L).$

\smallskip\noindent
\noindent
2. For $\mu= pd/2$ with an integer $p\geq q$, the support of the measure 
$\delta_\xi \circ_\mu\delta_\eta$ describes the set of possible singular spectra
of sums $x+y$ made up by matrices $x,y\in M_{p,q}(\b F)$ with given singular spectra
$\xi$ and $\eta$. 

\end{remarks}

\vspace{8pt}

Let us return to general indices $\mu\in \mathcal M_q.$
In analogy to Lemma \ref{dual} for orbit hypergroups from groups, we expect that the characters  of the hypergroup $\Xi_{q,\mu}$ 
are all obtained by taking  $U_q$-means of the characters of $\Pi_{q,\mu}$.
For $\xi \in \Xi_q$, define $\psi_\xi = \psi_\xi^\mu\in C_b(\Xi_q)$ by
\[\psi_\xi(\eta) :=\int_{U_q} \phi_\xi(u\eta u^{-1})du\]
where $\phi_\xi(r)= \phi_\xi^\mu(r) = \mathcal J_\mu(\frac{1}{4}r\xi^2r).$ The $\phi_\xi$  are the characters of $\Pi_{q,\mu}$ which are parametrized by diagonal matrices.
Note that by $U_q$-invariance of  $\mathcal J_\mu$ we have
\begin{equation}\label{charconj}
\phi_s(uru^{-1}) = \phi_{u^{-1}su}(r) \quad \forall\,r,s\in \Pi_q.
\end{equation}
Hence the mean of $\phi_s$ equals the mean of $\phi_\xi$ for $\xi=\sigma(s)$
and $\psi_\xi(\eta)= \psi_\eta(\xi)$ for
all $\xi, \eta\in \Xi_q$.

\begin{theorem}\label{dualchamber} Let $ \mu \in \mathcal M_q$.  
\begin{enumerate}
\item[\rm{(1)}] The dual space of the hypergroup $\,\Xi_{q,\mu} = (\Xi_q, \circ_\mu)$ is given by
\[ \widehat\Xi_{q,\mu}\,=\, \{\psi_\xi = \psi_\xi^\mu\,:  \xi \in \Xi_q\}.\]
\item[\rm{(2)}]
The hypergroup $\Xi_{q,\mu}$ is self-dual via the homeomorphism $\,\Xi_{q,\mu}\to \widehat\Xi_{q,\mu}, \, \xi\mapsto\psi_\xi$ . Under this identification, the Plancherel measure $\widetilde\pi_\mu$ of $\Xi_{q,\mu}$ coincides with the Haar measure $\widetilde\omega_\mu$. 
\end{enumerate}
\end{theorem}

\begin{proof}  Ad (1). It is easily checked that each $\psi_\xi$ is multiplicative w.r.t. $\circ_\mu$ and therefore belongs to  $\widehat\Xi_{q,\mu}$.
Indeed, for $\eta, \,\zeta\in \Xi_q$ we calculate
\[\psi_\xi(\eta\circ_\mu \zeta) = \int_{U_q} (\psi_\xi\circ\pi)(\eta*_\mu u\zeta u^{-1})du =\int_{\Pi_q}\int_{U_q}\int_{U_q} \phi_\xi(vrv^{-1})d(\delta_\eta *_\mu \delta_{u\zeta u^{-1}})(r)dudv.\]
By Lemma \ref{l1} this equals
\[ \int_{U_q}\int_{U_q} \phi_\xi(v\eta v^{-1} *_\mu vu\zeta u^{-1}v^{-1}) dudv \,=\, \psi_{\xi}(\eta) \psi_\xi(\zeta).\]
It remains to show that each character of $\Xi_{q,\mu}$ is of the form $\psi_\xi$ with some $\xi \in \Xi_q$. For this, notice first
that the hypergroup $\Xi_{q,\mu}$  has subexponential growth, just as $\Pi_{q,\mu}$. Thus according to Theorem 2.17 of \cite{V1},
the support of the Plancherel measure $\widetilde\pi_\mu$
of $\Xi_{q,\mu}$ coincides with the full dual $\widehat\Xi_{q,\mu}$. Let $\psi\in \widehat\Xi_{q,\mu} = \text{supp}\,\widetilde\pi_\mu$. Then by  Corollary 6 of \cite{V2} there exists a sequence of functions $f_n\in C_c(\Xi_q)$ such that $f_n \circ_\mu \overline f_n$ converges to $\psi$ locally uniformly. Hence
$(f_n\circ\pi)*_\mu(\overline f_n\circ\pi)$ converges to $\psi\circ\pi$  locally uniformly on $\Pi_q$, which implies that $\psi\circ\pi$ is positive definite 
on the hypergroup $\Pi_{q,\mu}$. 
Notice that $(\psi\circ\pi)(0) =1$. Thus 
by Bochner's Theorem for commutative hypergroups (Thm. 12.3.B of \cite{Je})
and the self-duality of $\Pi_{q,\mu}$, there exists a probability measure $a\in M^1(\Pi_q)$ such
that
\begin{equation}\label{charident}
(\psi\circ\pi)(r)\, =\, \int_{\Pi_q} \phi_s(r) da(s) \quad \forall \, r\in \Pi_q.
\end{equation}
As $\psi\circ\pi$ is invariant under the action of $U_q$ and in view of \eqref{charconj}, the measure $a$ must be $U_q$-invariant as well,
i.e. $\int f(uru^{-1})da(r) = \int f(r)da(r)$ for all $f\in C_c(\Pi_q)$. 
Let $\widetilde a$ denote the image measure of $a$ under $\pi$. Then with $\eta = \pi(r)$ relation \eqref{charident} becomes
\[ \psi(\eta) \,=\, \int_{\Pi_q}\phi_s(\eta)da(s) \,=\,\int_{\Xi_q} \psi_\xi(\eta)\, d\,\widetilde a(\xi) .\]
On the other hand, by Bochner's Theorem for the hypergroup $\Xi_{q,\mu}$, 
the character $\psi$ of $\Xi_{q,\mu}$ is an extremal point of the set
of positive definite functions $f$ on the hypergroup $\Xi_{q,\mu}$ with the
additional property  $f(0) =1$. This implies that
$\widetilde a$ must be a point measure, i.e. $\widetilde a = \delta_\xi$ for
some $\xi\in \Xi_q$. Hence $\psi = \psi_\xi$.

\medskip
Ad (2). The self-duality of $\Xi_{q,\mu}$ is proven in the same way as that of 
$\Pi_{q,\mu}$ (Theorem \ref{dualspace}). To determine the Plancherel measure, 
 let $f\in C_c(\Xi_q)$ and $g:= f\circ \pi$. Then the Fourier transforms
of $f$ and $g$ w.r.t. the hypergroup structures $\Xi_{q,\mu}$ and $\Pi_{q,\mu}$ are related via
\[ \widehat g(s) = \int_{\Pi_q} g(r) \phi_s(r) d\omega_\mu(r) \,=\, \int_{\Xi_q} f(\xi)
\psi_{\sigma(s)}(\xi)\,d\widetilde\omega_\mu(\xi)\, =\, (\widehat f\circ \pi)(s).\]
By the Plancherel theorem for $\Pi_{q,\mu}$ we readily obtain
$\, \int_{\Xi_q} |f|^2 d\widetilde\omega_\mu  \,=\, \int_{\Xi_q}  |\widehat f\,|^2 d\widetilde\omega_\mu.$
This finishes the proof.
\end{proof}

We would like to write the characters $\psi_\xi$ in a more explicit form.
Recall that
$\phi_s(r) = \mathcal J_\mu\bigl(\frac{1}{4}sr^2s\bigr)$ for $r,s\in \Pi_q$, 
where $\mathcal J_\mu$ is given in terms of the spherical series
\eqref{Besselseries}.
The spherical polynomials 
satisfy the product formula
\[ \frac{Z_\lambda(r)Z_\lambda(s)}{Z_\lambda(I)}  = \int_{U_q} Z_\lambda(\sqrt{r}usu^{-1}\sqrt{r})du \quad \forall \,r,s\in \Pi_q\,,\]
see \cite{FK}, Cor. XI.3.2. or \cite{GR1}, Prop. 5.5. 
 This 
implies an integral representation  for the Bessel functions $\mathcal J_\mu$ of two matrix arguments (recall Section \ref{BFM}):
\[\mathcal J_\mu(r,s) = \int_{U_q} \mathcal J_\mu(\sqrt{r}usu^{-1}\sqrt{r})du, \quad  r,s\in \Pi_q.\]
Thus for $\xi,\eta\in \Xi_q$ we have
\begin{equation}\label{charid2}\psi_\xi(\eta) = \,\int_{U_q} \mathcal J_\mu\bigl(\frac{1}{4}\xi u\eta^2u^{-1}\xi\bigr)du \,=\,
\mathcal J_\mu\bigl(\frac{\xi^2}{2}, \frac{\eta^2}{2}\bigr).
\end{equation}
We shall use this representation in order to identify the
characters $\psi_\xi$ with Dunkl type Bessel functions for the root system of type $B_q$.

\subsection{Bessel functions associated with root systems}\label{dunklbessel}

Bessel functions associated with root systems are an important ingredient in the 
theory of rational Dunkl operators, which was initiated by C.F. Dunkl in 
the late 80ies (\cite{D1}, \cite{D2}). They are a symmetrized version of the Dunkl kernel, which is the analogue of the usual exponential function in this theory. As a subclass, they
include the spherical 
functions of a Cartan motion group, c.f. \cite{dJ2}.
In this section we give a brief account on Dunkl theory and the associated Bessel functions; for a general background, the reader is referred to 
\cite{DX}, \cite{O} and \cite{R3}. 

Let $G$ be a finite reflection group  on $\b R^q$ (equipped with the usual
Euclidean scalar product $\langle \,.\,,\,.\,\rangle$), and let $R$ be the reduced root
system of $G$. We extend the action of $G$ to $\b C^q$ and 
$\langle \,.\,,\,.\,\rangle$ to a bilinear form on $\b C^q\times \b C^q$. 
A funciton $k: R\to \b C$ which is invariant under $G$
is called a multiplicity function on $R$. Important special cases of reflection groups are the symmetric group $S_q$ which acts on $\b R^q$ by permuting the standard basis vectors $e_i$, and the
hyperoctahedral group $B_q = S_q\ltimes \b Z_2^q$ which acts by permutations of the basis vectors and  sign changes. The root system of $B_q$ is given by $R =\{ e_i, \,1\leq i\leq q\} 
\cup \{ \pm e_i\pm e_j, 1\leq i<j\leq q\}$, and a multiplicity on it is of the form $k = (k_1, k_2)$ where $k_1$ is  the value on the 
roots $\pm e_i$ and $k_2$ is the value on
the roots $\pm e_i \pm e_j$.

For a finite reflection group $G$ and a fixed multiplicity function $k$ 
on its root system, the associated (rational)
Dunkl operators are defined by
\[ T_\xi = T_\xi(k) = \partial_{\xi} + \frac{1}{2}\sum_{\alpha\in R}
k_\alpha\langle\alpha,\xi\rangle \frac{1}{\langle\alpha,\,.\,\rangle}(1-\sigma_\alpha), \quad \xi \in \b C^q;\]
here $\sigma_\alpha$ denotes the reflection in the hyperplane perpendicular to $\alpha$ and the action of $G$ is extended to functions on $\b C^q$ via
$g.f(\xi)= f(g^{-1}\xi)$.
The  $T_\xi$ are  homogeneous of degree $-1$ on 
the space $\mathcal P =\b C[\b C^q]$ 
of polynomial functions on $\b C^q$ and they commute: $T_\xi T_\eta = T_\eta T_\xi$
(\cite{D1}). 
Hence the map $\xi\mapsto T_\xi$ extends uniquely to a linear map $p\mapsto p(T), \, \mathcal P\to \text{End}(\b C^q)$. 
The Dunkl operators induce a sesquilinear pairing 
\[ [p,q]_k = (p(T)\overline{q})(0)\]
on $\mathcal P$, where $\overline{q}(\xi) := \overline{q(\overline{\xi})}$.
In the following we assume that $k$ is non-negative. Then 
$[p,q]_k $ is actually a scalar product on $\mathcal P$, see Prop. 2.4. of \cite{DO}. Moreover, 
for each fixed $w\in \b C^q$, the joint eigenvalue problem 
\[T_\xi f\,=\, \langle\xi,w\rangle f \quad\forall\,\xi\in \b C^q; \quad f(0)=1\]
has a unique holomorphic solution $f(z) = E_k(z,w)$ called the 
Dunkl kernel. It is symmetric in its arguments and satisfies
$\, E_k(\lambda z,w) = E_k(z,\lambda w)$ for all $\lambda\in \b C$ as well as
$\, E_k(gz,w) = E_k(z,gw)$ for all $g\in G$. 
The generalized Bessel function
\[ J_k(z,w) := \frac{1}{|G|}\sum_{g\in G} E_k(z,gw)\]
is $G$-invariant in both arguments. Moreover, $g(z)= J_k(z,w)$ is the unique holomorphic solution of the "Bessel system"
\begin{equation}\label{besseleq}
p(T) g\, = \,p(w)g \quad\forall\, p\in \mathcal P^G; \quad g(0)=1
\end{equation}
where $\mathcal P^G$ denotes the subalgebra of $G$-invariant polynomials
in $\mathcal P$, see \cite{O}. For crystallographic reflection groups and certain values of $k$, the operators $p(T)$, when restricted to $G$-invariant functions on $\b R^q$, constitute  the system of invariant differential operators of a Euclidean-type symmetric space and the Bessel functions $J_k(\,.\,,w)$ can be identified with the associated spherical functions; for details see 
\cite{dJ2}. The Dunkl kernel $E_k$ gives rise to an integral transform on 
$\b R^q$ called the Dunkl transform. Let $w_k$ denote the  weight function 
\[ w_k(x) = \prod_{\alpha\in R} |\langle\alpha,x\rangle|^{2k_\alpha}\]
on $\b R^q$. The Dunkl transform is the integral transform on $L^1(\b R^q, w_k)$ defined by
\[ f \,\mapsto \widehat f^{\,k}, \qquad \widehat f^k(\xi) = c_k^{-1}\int_{\b R^q} f(x) E_k(-i\xi,x) w_k(x) dx \quad (\xi\in \b R^q)\]
with the constant
\[ c_k := \int_{\b R^q} e^{-|x|^2/2} w_k(x) dx.\]
A thorough study of this transform is given in \cite{dJ1}. It has many properties in common with the usual Fourier transform to which it reduces in case
$k=0$.
In particular, the Dunkl transform (as normalized above) extends to an isometric isometric isomorphism of $L^2(\b R^q, w_k)$, and $\,(T_\eta f)^{\wedge k}(\xi) =
i\langle \xi,\eta\rangle\widehat f^{\,k}(\xi)$ 
for differentiable $f$ of sufficient decay. 
It is a long-standing open question whether  $L^1(\b R^q, w_k)$ 
can be given the structure of a commutative Banach algebra so that the Dunkl transform becomes the Gelfand transform on its (symmetric) spectrum,
similar as for commutative hypergroups. In the rank-one case there is such
a convolution, but it is not positivity-preserving. For details  and affirmative results in this direction see \cite{R2}. It is however conjectured that for arbitrary $G$ and $k\geq 0$, the Bessel functions $J_k$ have a positive product formula
which leads to a commutative hypergroup structure on a distinguished closed 
Weyl chamber $\Xi$ of $G$, the dual of this hypergroup being made up by the functions $\xi\mapsto J_k(\xi, \eta), \, \eta\in \Xi$. In rank one and  in all Cartan motion group cases this is true, see \cite{R2}.
In the following, we shall confirm this conjecture for three
continuous series of multiplicities for root system $B_q$.
Indeed, we shall identify the characters of the hypergroups $\Xi_{q,\mu}\,$ in Section \ref{chamberconvo} with 
Dunkl-type  Bessel functions for $B_q$ and thus obtain
hypergroup structures with these Bessel functions as characters.

\subsection{Dunkl theory and the convolutions on the Weyl chamber}
\label{connect}

Let us denote by $J_k^B$  the Dunkl-type Bessel function associated
with the reflection group $G=B_q$ and multiplicity $k = (k_1,k_2)$, and by $[\,.\,,\,.\,]_k^B$ the associated Dunkl pairing. For $z=(z_1, \ldots, z_q)\in \b C^q$ we put $z^2= (z_1^2, \ldots, z_q^2)$. The following key result identifies $J_k^B$ with a generalized $_0F_1$-hypergeometric function of two arguments
(recall the notions of Section \ref{BFM}):

\begin{proposition} Let $k=(k_1, k_2)\geq 0$ and $k_2 >0$. Then
for all $z,w\in \b C^q$, 
\[ J_k^B(z,w) = \,_0F_1^\alpha\bigl(\mu;\frac{z^2}{2},\frac{w^2}{2}\bigr)
\quad\text{ with }\,\,\alpha = \frac{1}{k_2}\,, \,\mu = k_1 +(m-1)k_2 + \frac{1}{2}.\]
\end{proposition}

This result was already 
mentioned in Section 6 of \cite{BF}, but the reasoning there is rather sketchy, and there is an erraneous sign in one of the 
arguments. We therefore include a  proof by different methods.

\begin{proof} The
modified Jack polynomials $p_\lambda(z) = C_\lambda^\alpha(z^2),$ indexed by partitions $\lambda\geq 0$,
are homogeneous of degree $2|\lambda|$ and form a basis of the 
vector space $\mathcal P^G$ for $G= B_q$. 
Thus the  Bessel function $J_k^B$ has a homogeneous expansion 
of the form 
\begin{equation}\label{Besselexpan} J_k^B(z,w)\,=\,\sum_{\lambda\geq 0} a_\lambda(w)p_\lambda(z)
\end{equation}
with certain coefficients $a_\lambda(w)\in \b C$. 
In view of the Bessel system we have
\[ p_\lambda(T^z) J_k^B(z,w)\big\vert_{z=0} \, =\, p_\lambda(w) J_k^B(0,w) \,=\, C_\lambda^\alpha (w^2),\]
where the superscript $z$ indicates operation w.r.t. the variable $z$.
On the other hand, the results of \cite{BF} (relation (2.9) and the 
formula on top of p.~214) imply that the $p_\lambda$ are orthogonal with respect to  $[\,.\,,\,.\,]_k^B$ with 
\[ [p_\lambda, p_\lambda]_k^B \,=\, 4^{|\lambda|} |\lambda|!\, (\mu)_\lambda^\alpha \cdot C_\lambda^\alpha({\bf 1})\,=:M_\lambda,  \quad \mu = k_1 + (m-1) k_2 + \frac{1}{2}\,.\]
Differentiation of \eqref{Besselexpan} (recall that $p_\nu$ has real coefficients) now gives
\[ p_\lambda(T^z) J_k^B(z,w)\big\vert_{z=0} \,=\, \sum_{\nu\geq 0} a_\lambda(w)
[p_\lambda\,,p_\nu]_k^B \,=\, M_\lambda a_\lambda(w).\]
Hence $\,a_\lambda(w) = M_\lambda^{-1}\,C_\lambda^\alpha(w^2),$ which implies the assertion.
\end{proof}

As a consequence of \eqref{dunklbessel2},  Bessel functions associated with a symmetric cone
can now be identified with Dunkl Bessel functions of type $B_q$ with specific multiplicities:

\begin{corollary} \label{cor1} Let $\Omega$ be an irreducible symmetric cone inside a  Euclidean 
Jordan algebra of rank $q$. Then 
 for $r,s\in \overline{\Omega}$ with eigenvalues $\,\xi=( \xi_1, \ldots ,\xi_q)$ and 
$\eta=(\eta_1, \ldots,\eta_q)$ respectively, we have
\[ \mathcal J_\mu\bigl(\frac{r^2}{2},\frac{s^2}{2}\bigr) = J_k^B(\xi,i\eta)\]
where $k$ is given by $\,k = k(\mu,d) = \bigl(\mu-\frac{d}{2}(q-1)-\frac{1}{2}, \, \frac{d}{2}\bigr).$
\end{corollary}

\noindent
Now consider again the hypergroup structures $\Xi_{q,\mu} = (\Xi_q, \circ_\mu)$ on 
\[\Xi_q = \{\xi =(\xi_1,\ldots,\xi_q)\in \b R^q: \xi_1\geq \ldots\geq \xi_q\geq 0\}\] which is actually a closed Weyl chamber for the reflection group $B_q$. 
The consequence of our above identification can be formulated in a twofold way: 

\begin{corollary} \label{cor2}   The characters of the hypergroup $\Xi_{q,\mu}$,  $\mu\in \mathcal M_q$ 
 are given by 
\[ \psi_\eta(\xi) = J_k^B(\xi,i\eta), \quad \eta\in \Xi_q,\]
with the multiplicity $k=k(\mu,d)$ as in the previous corollary. 
\end{corollary}

\begin{corollary}\label{cor3}
Consider the root system of type $B_q$ with
a  multiplicity $k= (k_1, k_2)$ where  $k_2 = \frac{d}{2}$ with $d\in \{1,2,4\}$ and $k_1= \frac{d}{2}(p-q+1) -\frac{1}{2}$ for integer $p\geq q$ or arbitrary $k_1 \geq \frac{1}{2}(dq-1)$. Then the associated 
Dunkl-type Bessel functions 
$\xi\mapsto J_k^B(\xi,i\eta)$ are the characters of the  hypergroup
$(\Xi_q, \circ_\mu)$ 
on the closed Weyl chamber $\Xi_q$, where $\mu =  k_1 + (q-1)k_2 + \frac{1}{2}$ and the convolution $\circ_\mu$ is defined over $\b F = \b R, \b C, \b H,$ 
depending on the value of $d$. In particular, 
 the
Bessel function $J_k^B$ satisfies the positive product formula
\[ J_k^B(\xi,z) J_k^B(\eta,z) \,=\, \int_{\Xi_q} J_k^B(\zeta, z) \,d(\delta_\xi \circ_\mu\delta_\eta)(\zeta) \quad \forall \, \xi,\eta\in \Xi_q, 
\,z\in \b C^q.\]
\end{corollary}  

The
 hypergroup Fourier transform on $\Xi_{q,\mu}$ is  given by
\[ \widehat f(\eta) \,=\, \int_{\Xi_q} f(\xi) J_k^B(\xi,i\eta)\, d\widetilde\omega_\mu(\xi), \]
with $\widetilde\omega_\mu = d_\mu h_\mu(\xi)d\xi\,$ as in Theorem \ref{main3}.
Notice that $h_\mu$ coincides up to a constant factor with the weight
$w_k$ for  $k=k(\mu,d)$. As
 $w_k$ is $B_q$-invariant, we therefore have
 \[\widehat f(\eta) = \, const\cdot\widehat  F^{\,k}(\eta),\]
 where $F$ denotes
the $B_q$-invariant extension of  $f$ to $\b R^q$ and  $\widehat  F^{\,k}$ its Dunkl transform. Notice that
$\widehat  F^{\,k}$ is $B_q$-invariant as well.
With the Plancherel Theorem for the Dunkl transform at hand, we are now also in a position to determine
the 
normalization constant $d_\mu$ of $\widetilde\omega_\mu$ as announced previously.
Indeed, recall from Theorem \ref{dualchamber} that the Plancherel measure of the hypergroup $\Xi_{q,\mu}$ coincides with $\widetilde\omega_\mu$ under the natural identification of $\Xi_{q,\mu}$ with its dual. Using this and the Plancherel theorem for the Dunkl transform,   we readily obtain 
\[ \widehat f \,=\, \widehat F^{\,k} \vert_{\,\Xi_q} \]
and 
\[ d_\mu = \Bigl(\int_{\Xi_q} h_\mu(x) e^{-|x|^2/2}dx\Bigr)^{-1}.\]
The value of $d_\mu$ can be calculated explicitly; it is a particular case of a Selberg type integral which was evaluated by Macdonald 
\cite{M1} for the classical root systems.

In the general Dunkl setting, there is a generalized translation
on suitable function spaces  which replaces
the usual group addition to some extent, see \cite{R1}, \cite{R2} and the references cited there. On $L^2(\b R^q, w_k)$, this translation is defined by
\[ \tau_\eta f(\xi) =\, c_k^{-1} \int_{\b R^q} \widehat f^{\,k}(\xi)E_k(i\xi,\zeta) E_k(i\eta,\zeta) w_k(\zeta) d\zeta.\]
On has $\tau_\eta: L^2(\b R^q, w_k)\to L^2(\b R^q, w_k)$ with $\, (\tau_\eta f)^{\wedge k}(\zeta) = E_k(i\eta, \zeta) \widehat f^{\,k} (\zeta).$ 
If we restrict to the  Weyl group invariant case for $B_q$ with multiplicities as in Corollary \ref{cor3}, then  this generalized translation just coincides with the translation defined in terms of hypergroup convolution. If, say, $f$ belongs to $L^2(\b R^q, w_k)$ and is also continuous and Weyl group invariant, 
then we have with the notions of Corollary \ref{cor3}
\[ \tau_\eta f(\xi) \,=\, \delta_\xi *_\mu \delta_\eta(f) \quad \forall \xi,\eta \in \Xi_q. \]

When $\mu = \frac{pd}{2}$ with an integer $p\geq q$, i.e. $\,k_1= \frac{d}{2}(p-q+1) -\frac{1}{2}$, then the Bessel functions $J_k^B(\,.\,,z)$  can be identified with the 
the spherical functions of the Cartan motion group associated with the Grassmann manifold
$ U(p,q)/(U_p\times U_q)$. 
This follows from the discussion in \cite{dJ2} (see also \cite{R2}),
and is in accordance with Remark \ref{double} which implies that for  $\mu=\frac{pd}{2}$, 
the hypergroup convolution of $\Xi_{q,\mu}$ coincides with that of biinvariant 
measures for the Gelfand pair $((U_p\times U_q)\ltimes M_{p,q},U_p\times U_q)$. 
The multiplicative functions coincide with  the (elementary) spherical functions of this Gelfand pair. Thus the hypergroups $\Xi_{q,\mu}$ with $\mu\geq \rho -1$ 
interpolate the discrete series of convolution algebras derived from the tangent space analysis on 
Grassmann manifolds.

\medskip\noindent
\emph{Acknowledgement.} The paper was for the most part written while the author was staying at the Korteweg-de Vries Institute for Mathematics at 
 the University van Amsterdam, being supported by the Netherlands Organisation for Scientific Research (NWO). It is a pleasure to thank the Korteweg-de Vries Institute at 
 the University van Amsterdam for their kind hospitality. Also, the support of 
the SFB/TR12  at the University of Bochum is greatly acknowledged.
 Special thanks go to Tom Koornwinder and to Michael Voit for various helpful comments and 
discussions.

\end{document}